\documentstyle[amsmath,amssymb]{article}

\def\be{\begin{eqnarray}}
\def\ee{\end{eqnarray}}

\hoffset= - 1.0in         % switch off for draft style
\begin{document}

\hfill ITEP/TH-13/07

\bigskip

\centerline{\Large{$A_\infty$ structure on simplicial complexes
}}

\bigskip

\centerline{\it V.Dolotin, A.Morozov and Sh.Shakirov}

\bigskip

\centerline{ITEP, Moscow, Russia}

\bigskip

\centerline{ABSTRACT}

\bigskip

A discrete (finite-difference) analogue of differential forms
is considered, defined on simplicial complexes, including
triangulations of continuous manifolds.
Various operations are explicitly defined on these forms,
including exterior derivative and exterior product.
The latter one is non-associative. Instead, as anticipated,
it is a part of non-trivial $A_\infty$ structure,
involving a chain of poly-linear operations, constrained by
nilpotency relation:
$\big(d + \wedge + m + \ldots\big)^n = 0$ with $n=2$.

\bigskip

\bigskip

\section{Introduction}

In the space of ordinary differential forms two operations
are naturally defined: linear exterior derivative $d$ 
and bilinear exterior multiplication $\wedge$.
They satisfy a triple of quadratic relations:
$$
\begin{array}{lclc}
d^2=0 && {\rm i.e.\ nilpotency\ of}\ d: & d(df)=0
\ \ \forall\ {\rm form}\ f\\
d\wedge + \wedge d = 0 && {\rm Leibnitz\ rule}:
%\ for\ external\ multiplication}:
& d(f\wedge g) = df \wedge g + (-1)^{|f|}f\wedge dg
\ \ \forall\ f,g \\
\wedge^2 = 0 && {\rm associativity\ of}\ \wedge: &
(f\wedge g)\wedge h = f\wedge (g\wedge h)
\ \ \forall\ f,g,h
\end{array}
$$
which can be rewritten as a single nilpotency constraint
\be
\big( d + \wedge \big)^2 = 0
\label{nccont}
\ee
Eq.(\ref{nccont}) is often used as a basic relation in
construction of topological field theories \cite{top}.
A natural {\it deformation} of (\ref{nccont}) in
BRST-BV formalism \cite{BV} is the $A_\infty$-structure
\cite{A8}, which can be symbolically written as
\be
\big( d + \wedge + m + \ldots \big)^2 =
\big(m^{(1)} + m^{(2)} + m^{(3)} + \ldots\big)^2 = 0
\label{A8}
\ee
and involves an infinite chain of $k$-linear operations
with all natural $k$.
Exterior derivative
$d=m^{(1)}$ and exterior product $\wedge = m^{(2)}$
are just the first two members of this chain.
More accurately, (\ref{A8}) is an infinite collection of
quadratic constraints:
\be
\sum_{k + l = n + 1} m^{(k)} m^{(l)} = 0 \ \ \ \ n = 1, 2, 3, \ldots
\label{RMK}
\ee
i.e.
\be
d^2 = 0
\label{A2}
\ee
\be
d\wedge + \wedge d = 0
\label{A3}
\ee
\be
dm + md + \wedge^2 = 0
\label{A4}
\ee

\begin{center}
\ldots
\end{center}
It is a long-standing claim \cite{Lo} that the
$A_\infty$-structure is a natural property of
discrete de Rham complex, i.e. arises in the study of
finite-difference forms on simplicial complexes \cite{Pont},
in particular, on triangulations of continuous manifolds.
Despite certain efforts and progress in investigation
of this idea (\cite{prog1}, \cite{prog2}) it does not attract attention
that it deserves and explicit realization is still lacking.
Once established, such construction
should possess further generalization
to higher-order nilpotents, in the spirit of
non-linear algebra \cite{nolal}.
Such generalizations are expected to provide new tools
for a variety of applications, from topology to strings.

It is the goal of this paper to present explicit formulas
about the forms on simplicial complexes,
operations on the space of these forms and relations between
these operations.
Such down-to-earth description of the lowest components of
the $A_\infty$-structure shows what {\it explicit} realization
means and can hopefully stimulate the search
for equally explicit formulas for its cubic and higher-order
generalizations.
\pagebreak

\section{Basic example in dimension one
\label{BE}}
We begin with a one-dimensional example of discrete forms
and operations. Consider a one-dimensional discrete space
(a graph).
Such a space contains points (vertices) and links between
some of them. Function $f_i$ is defined on its vertices $i$,
while 1-form $\psi_{ij}$ is defined on links $ij$.

\subsection{Differential $d$}

Discrete derivative of a function $df$ is a 1-form:
for a given link $ij$, connecting the points $i$ and $j$,
it is a difference between the values of $f$ at the two ends
of a link:
\be
df_{ij} = f_j - f_i
\ee
In one dimension, differential of a 1-form is zero,
and so the $d^{2} = 0$ identity is trivial.

\subsection{Product $\wedge$}

The product of two functions is the point-wise product:
\be
(f \wedge g)_i = f_ig_i\ \ \ \ \forall\ i
\ee
A natural definition for a product of a function $f$ and
a $1$-form $\psi$ is again an $1$-form,
\be
(f \wedge \psi)_{ij} = (\psi \wedge f)_{ij} =
\frac{1}{2}(f_i+f_j)\psi_{ij}\ \ \ \ \forall\ ij
\ee
In one dimension, a product of two 1-forms is zero.

\subsection{The Leibnitz rule $d \wedge + \wedge d = 0$}
Operations $d$ and $\wedge$ satisfy the usual Leibnitz rule
\be
d(f \wedge g) = df \wedge g + f \wedge dg
\label{deta}
\ee
The calculation is straightforward (f and g are functions):
$$\big(d(f \wedge g)\big)_{ij} = (f \wedge g)_j - (f \wedge g)_i
= f_jg_j - f_ig_i,$$

$$\big(df \wedge g\big)_{ij} = \frac{1}{2}(df)_{ij}(g_i+g_j) =
\frac{1}{2}(f_j-f_i)(g_i+g_j)$$

$$\big(f \wedge dg\big)_{ij} = \frac{1}{2}(f_i+f_j)(dg)_{ij} =
\frac{1}{2}(f_i+f_j)(g_j-g_i)$$

$$\big(f \wedge dg\big)_{ij} + \big(df \wedge g\big)_{ij} =
\frac{1}{2}(f_i+f_j)(g_j-g_i) + \frac{1}{2}(f_j-f_i)(g_i+g_j) =
f_jg_j - f_ig_i = \big(d(f \wedge g)\big)_{ij}$$
If we define differential not only on functions,
but also on pairs of functions as
\be
d(f,g) = - (df,g) - (f,dg)
\ee
then the Leibnitz rule (\ref{deta}) can be represented in three
different ways:
$$d(f \wedge g) - df \wedge g - f \wedge dg = 0$$
$$d\big{(}\wedge (f,g)\big{)} - \wedge(df, g + f,dg) = 0$$
$$(d \circ \wedge + \wedge \circ d) (f,g) = 0$$
and, finally
\be
d \wedge + \wedge d = 0
\ee

In this paper, we often use this kind of notation.
We call this form of $A_{\infty}$-relations \textit{brief},
in contrast with the \textit{detailed} form (\ref{deta}).
The brief form, implied by the BRST-BV formalism,
is extremely useful,
both from practical and theoretical points of view.
First of all, higher $A_{\infty}$-relations become
more and more complicated, and one needs the brief form
to handle them in a simple and transparent way.
Second, in the {\it brief} form the $A_\infty$ structure
admits a simple interpretation as deformation of $d$
by $\wedge$, see s.\ref{dco} below.

\subsection{Non-associativity $\wedge^{2} \neq 0$}

The discrete product is not obligatory associative. While
$(f \wedge \psi) \wedge g = f \wedge (\psi \wedge g)$,
for the same triple in another order
$(f \wedge g) \wedge \psi \neq f \wedge(g \wedge \psi)$.
Indeed,
$$
\Big((f\wedge\psi\big)\wedge g\Big)_{ij} -
\Big(f\wedge (\psi\wedge g)\Big)_{ij} = $$
\be
\frac{1}{2}(f\wedge\psi)_{ij}(g_i+g_j) -
\frac{1}{2}(f_i+f_j)(\psi\wedge g)_{ij} =
\frac{1}{4}(f_i+f_j)\psi_{ij}(g_i+g_j)
- \frac{1}{4}(f_i+f_j)\psi_{ij}(g_i+g_j) = 0,
\ee
but
$$
\Big(\big(f\wedge g)\wedge\psi\Big)_{ij} -
\Big(f\wedge(g\wedge\psi)\Big)_{ij} =
$$ $$
\frac{1}{2}\big((f\wedge g)_i+(f\wedge g)_j\big)\psi_{ij} -
\frac{1}{2}(f_i+f_j)(g\wedge\psi)_{ij} =
\frac{1}{2}(f_ig_i+f_jg_j)\psi_{ij}
-\frac{1}{4}(f_i+f_j)(g_i+g_j)\psi_{ij} =
$$
\be
=\frac{1}{4}(f_i-f_j)(g_i-g_j)\psi_{ij}
\ee
is not vanishing unless $df = 0$ or $dg = 0$.
This is a fact of great importance:
in general, a discrete product satisfying the Leibnitz rule
fails to be associative.

If we define the product $\wedge$ not only on pairs of forms,
but also on triples (converting them into pairs):
\be
\wedge \big(a,b,c\big) =
\big(a,\wedge(b,c)\big) - \big(\wedge(a,b),c\big)
\ee
then
$$
\big{(} a \wedge (b \wedge c) - (a \wedge b) \wedge c \big{)} =
\wedge \big{(} a,\wedge(b,c) \big{)} - \wedge \big{(}\wedge(a,b),c\big{)} = $$
\be
= \wedge \big{(} a,\wedge(b,c) - \wedge(a,b),c \big{)} =
\wedge \circ \wedge (a,b,c) = \wedge^{2} (a,b,c)
\ee
In this notation non-associativity simply means that
$\wedge^{2} \neq 0$.
This is just another example of the brief form of $A_\infty$
relations:
simplification occurs when we define operations
on a larger set of arguments.

\subsection{Higher operation $m$}

Consider a $3$-linear operation $m$ which converts a triple of
a function $f$  and two $1$-forms $\psi$ and $\chi$
into a $1$-form
\be
m(f,\psi,\chi)_{ij} = m(f,\chi,\psi)
= \frac{1}{8}(f_j-f_i)\psi_{ij}\chi_{ij}
\label{dim011}
\ee
The following (non)associativity relations hold:
$$
(f\wedge\psi)\wedge g - f\wedge(\psi\wedge g)
= m(df,\psi,g) - m(f,\psi,dg) = 0,
$$
\be
(f\wedge g)\wedge\psi - f\wedge(g\wedge\psi)
= m(df,g,\psi) + m(f,dg,\psi) \neq 0.
\label{mop}
\ee
Complemented by $m(f,g,\psi)_{i} = m(f,\psi,g)_i = 0$,
they provide examples of the first non-trivial relation (\ref{A4})
of the $A_\infty$ structure.
This relation states that for any forms $a,b,c$
\be
\Big(d \circ m + m \circ d + \wedge \circ \wedge\Big)
(a, b, c)= 0
\label{dmmd}
\ee
where
$d(a,b,c) = (da,b,c) + (-1)^{|a|} \cdot (a,db,c)
+ (-1)^{|a| + |b|}\cdot (a,b,dc)$.
Degree of a function $|f| = 0$, while for a 1-form $|\psi| = 1$.
In full detail eq.(\ref{dmmd}) states that
\be
dm(a,b,c) + m(da,b,c) + (-1)^{|a|} m(a,db,c) +
(-1)^{|a| + |b|} m(a,b,dc) = \big{(}(a \wedge b) \wedge c - a \wedge (b \wedge c)\big{)}
\label{m3}
\ee
This example is simple: in one dimension only 0-forms and 1-forms
are present.
However, it is informative: already at this level discrete
non-associativity can be observed,
together with the higher $A_{\infty}$ operation
$m = m^{(3)}$.
Non-associativity is not an isolated phenomenon,
it is an element of non-trivial $A_{\infty}$ structure
of discrete algebra of operations.
Non-associativity is fully described by the operation $m^{(3)}$:
instead of $\wedge^{2} = 0$ a more complicated identity
$$\wedge^2 = - d m - m d$$ holds.
The rest of the paper is devoted to generalization of this example
to arbitrary forms on arbitrary simplicial complexes.

\pagebreak
\section{Logic and plan of the presentation}

We proceed in four steps.

\subsection{Simplicial complexes, forms and operations}

{\bf First,} in section 4 we introduce the three "starting"
operations $m^{(1)}=d$,
its conjugate $\partial = d^\dagger$ and
$m^{(2)} = \wedge$, \linebreak which act on the space
$\Omega_*(M)$ of {\it forms}
on an {\it ordered simplicial complex} $M$.
Without ordering, $M$ can be considered as
equilateral triangulation of a manifold
(perhaps, singular) or simply as triangular
lattice in arbitrary dimension.
Some of the definitions in s.4
will be given in a form, which allows and suggests
further generalizations, especially to
higher nilpotent operators $d_n$:
\be
d_n^n = 0
\ee
In the main part of this paper we consider
only $d = d_2$, details of the $n>2$ case will be
presented elsewhere \cite{hini}.

\subsection{Formulation of the problem}

{\bf Second,} in section 5 we turn to the $A_\infty$ structure.
In this paper we do not discuss its origins and
hidden meanings, even cohomological and BRST-BV
aspects of the story will not be considered.
Instead our goal is to {\it construct} operations
$m^{(p)}$ on simplicial complex in the most explicit way,
and (\ref{A8}) is used as {\it defining property}
of these operations, {\it not} as a theorem derived for
some given "natural" operations.
In order to use (\ref{RMK}) in this quality, we
rewrite it as a sequence of recurrent relations:
\be
\boxed{\ \
dm^{(p)} + m^{(p)}d =
{\cal M}^{(p)} \equiv
- \sum_{q=2}^{p-1} m^{(q)} m^{(p+1-q)},\ \ \ p\geq 1
\ }
\label{A8''}
\ee
i.e.
${\cal M}^{(1)} = {\cal M}^{(2)} = 0,\ $
${\cal M}^{(3)} = - \wedge^2,\ $
${\cal M}^{(4)} = - \wedge m^{(3)} - m^{(3)}\wedge \ $
and so on.

Now the problem is formulated,
and the first thing to do is to check its consistency:
for (\ref{A8''}) to have a solution $m^{(p)}$ the r.h.s.
should satisfy
\be
d {\cal M}^{(p)} = {\cal M}^{(p)} d
\label{cons}
\ee
This consistency condition is considered in section \ref{concon}.

If (\ref{cons}) is true, solution to (\ref{A8''}) can
exist, but it does not need to be unique -- and it
actually {\it is} ambiguous. We consider two solutions: one naive,
but somewhat
distractive, in s.6, another nice and
fully satisfactory, but somewhat sophisticated,
in s.7. These will be our steps {\bf three}
and {\bf four}.

\subsection{Naive solution with the help of $K$-operator}

At {\bf step three}, in section 6,
we resolve eq.(\ref{A8''}) with
the help of the $d^{-1}$-operator $K$ \cite{KOP},
\be
dK + Kd = I
\label{Kop}
\ee
Such operator exists if $d$ is cohomologically trivial.
Even for topologically trivial simplicial complexes
there could be a problem with this at the level of $0$-forms
(since zeroth homology does not vanish) and we use a special trick
-- postulate
\be
\partial ({\rm point}) = \varnothing \neq 0
\label{postu}
\ee
-- to cure this problem and make
all formulas work universally on entire $\Omega_*$.
$K$ operator can be constructed from arbitrary operator
$\alpha$, provided associated "Laplace operator"
$\Delta_\alpha \equiv d\alpha + \alpha d$ is invertible:
\be
K_\alpha = \frac{1}{\Delta_\alpha}\,\alpha
\label{Kalp}
\ee
In most applications one take $\alpha = \partial$, so that
$\Delta_\partial = dd^+ + d^+ d\ $ or simple $\Delta$
is the ordinary Laplace
operator.

Given $K_\partial$, one can immediately resolve (\ref{A8''}):
\be
\boxed{\
m^{(\!p\,)} = K_\partial {\cal M}^{(\!p\,)} =
\frac{1}{\Delta}\,\partial {\cal M}^{(\!p\,)}
}
\label{m=KM}
\ee
and this is our "naive" solution.
It is indeed a solution:
\be
dm^{(\!p\,)} \ \stackrel{(\ref{m=KM})}{=}\
dK{\cal M}^{(\!p\,)} \ \stackrel{(\ref{Kop})}{=}\
{\cal M}^{(\!p\,)} - Kd{\cal M}^{(\!p\,)}
\ \stackrel{(\ref{cons})}{=}\
{\cal M}^{(\!p\,)} - K{\cal M}^{(\!p\,)} d
\ \stackrel{(\ref{m=KM})}{=}\
{\cal M}^{(\!p\,)} - m^{(\!p\,)}d,
\label{cheso}
\ee
as requested by (\ref{A8''}).

For a given simplicial complex $M$ eq.(\ref{m=KM})
is an absolutely explicit formula.
Moreover for
a {\it finite} complex with finite number $\#(M)$ of points
any operator is represented by a finite-size matrix,
and therefore satisfies a finite-degree algebraic
equation.
It follows that inverse $\Delta^{-1}$ is actually a
finite-degree polynomial of $\Delta$. Eq.(\ref{m=KM}) provides
quite reasonable formulas for a complex $M$ which consists of a
single simplex,
they are especially nice because on a simplex
the Laplace operator trivializes,
\be
\left.\Delta\right|_{{\rm simplex}} =
{\rm multiplication\ by\ }\#({\rm points})
\label{Lasi}
\ee
However, in general $\Delta$ is given by a
highly non-local expression
that depends in sophisticated way on simplicial complex
${M}$: a minor modification of $M$, say,
addition of an extra simplex at the boundary, changes
all the formulas. Thus, despite (\ref{m=KM})
formally resolves our problem,
this is not the desired solution; we want another one,
which is {\it local}: does not depend on the structure
of $M$ far away from the place where
$m^{(\!p\,)}$ is evaluated.
Such local solution exists, but its construction
is somewhat involved and includes a number of new
notions and ideas.

\subsection{Local solution and the local K-operator}

This construction of {\it local} operations $m^{(\!p\,)}$
in section 7 is our {\bf step four}.

The first idea of this construction is to change (\ref{m=KM})
for
\be
m^{(\!p\,)} =  {\cal M}^{(\!p\,)}K_\partial
\label{m=MK}
\ee
The difference is that $K$-operator in (\ref{m=KM})
is acting on forms, i.e. on the space $\Omega_*$,
while $K$-operator in (\ref{m=MK}) acts in the space
$\Omega_*^{\otimes p}$, i.e. on the sets of several,
namely of $p$, forms, that is on $p$-uples of forms.

At this stage we need the notion of
{\it lifting}, i.e, defining the action of
$d$, $\partial$, $\wedge$, $\Delta$ and all
operations $m^{(p)}$ on higher tensor spaces $\Omega_*^{\otimes p}$
once it is defined on $\Omega_*$.
For example,
$$
d\big(\omega_1,\omega_2,\ldots,\omega_p\big) \sim
\big(d\omega_1,\omega_2,\ldots,\omega_p\big) +
(-1)^{|\omega_1|}\big(\omega_1,d\omega_2,\ldots,\omega_p\big)
+ \ldots + (-1)^{|\omega_1|+\ldots+|\omega_{p-1}|}
\big(\omega_1,\omega_2,\ldots,d\omega_p\big),
$$
$$
\partial\big(\omega_1,\omega_2,\ldots,\omega_p\big) \sim
\big(\partial\omega_1,\omega_2,\ldots,\omega_p\big) +
(-1)^{|\omega_1|}
\big(\omega_1,\partial\omega_2,\ldots,\omega_p\big)
+ \ldots + (-1)^{|\omega_1|+\ldots+|\omega_{p-1}|}
\big(\omega_1,\omega_2,\ldots,\partial\omega_p\big),
$$
and even
\be
\Delta\big(\omega_1,\omega_2,\ldots,\omega_p\big) =
\big(\Delta\omega_1,\omega_2,\ldots,\omega_p\big) +
\big(\omega_1,\Delta\omega_2,\ldots,\omega_p\big)
+ \big(\omega_1,\omega_2,\ldots,\Delta\omega_p\big).
\label{comuLa}
\ee
However, (\ref{m=MK}) is not much better than
(\ref{m=KM}), because it is still not local.

The real advantage of (\ref{m=MK}) is that it has
straightforward local modifications.
On the space $\Omega_*^{\otimes p}$ with $p\geq 2$ one acquires
a new possibility: one can {\it localize} operation
$\partial$, which was {\it totally non-local} on
$\Omega_*$.
In section 7.1 we introduce this localized
$[\partial]$, its conjugate $[d]$ and associated
modifications of Laplace operator
\be
\Delta_{loc} = [d][\partial] + [\partial][d]
\label{loLa}
\ee
and $K$-operator $[K] = [\partial] \Delta^{-1}_{loc}$.
Then
\be
\boxed{\
m^{(\!p\,)} = {\cal M}^{(\!p\,)} [K] =
{\cal M}^{(\!p\,)} \frac{[\partial]}{\Delta_{loc}}
\ }
\label{m=M[K]}
\ee
acting from $\Omega_*^{\otimes p}(M)$ to $\Omega_*(M)$,
is our desired local expression for $m^{(p)}$.
In section 7.5 we justify the validity of this formula.

In addition to being local, eq.(\ref{m=M[K]}) appears
to be {\it simple} and {\it nice}.
This is because of distinguished properties of
the localized Laplace operator $\Delta_{loc}$,
considered in section 7.2.
This operator is actually acting within particular
simplexes, $\Delta_{loc} = \left.\Delta_{loc}\right|_{{\rm simplex}}$,
but it does not coincide with
$\left.\Delta\right|_{{\rm simplex}}$ and does not
satisfy (\ref{Lasi}) and even (\ref{comuLa}).
Still, it is a remarkably simple operator:
its eigenvalues, though not all coincident as for
(\ref{Lasi}), are subsequent integers and inverse operator
is easy to construct.
This makes (\ref{m=M[K]}) a very explicit formula,
for every particular complex it can be effectively
handled by MAPLE or Mathematica.
It goes without saying that eqs.(\ref{dim011}) and
(\ref{mop}) are immediately reproduced by
(\ref{m=M[K]}).
Higher operations are described in section 7.7.

Moreover, as explained in s.7.7,
the entire $A_\infty$ structure can be considered as a
$\wedge$-induced conjugation of the nilpotent operator $d$
-- and this interpretation survives generalization to
higher nilpotents $d_n$.

Finally, in the Appendix section 9 we present a
number of explicit calculations with discrete forms.
These include examples concerning discrete differential,
discrete wedge multiplication, higher operations,
Laplace operators and more.

\section{Definitions}

\subsection{Discrete space}

A \textit{simplicial complex}
is a finite collection of sets
$M = \{\sigma_{1},\sigma_{2},...\}$,
where $\sigma \subset \sigma'$ and $\sigma' \in M$
implies $\sigma \in M$.
\linebreak The elementary sets $\sigma$ are called
\textit{simplexes}.
Complex is thought to be a ``discrete space'' --
a number of vertices ($0$-simplexes),
with some pairs connected by links ($1$-simplexes),
some triples -- by triangles ($2$-simplexes),
some quadruples -- by tetrahedra ($3$-simplexes) and so on.
%\linebreak
Of course, the simplest thing to do is to connect everything:
consider a set of $N \geq 1$ points together with {\it all}
its subsets.
This is a closed simplex.
For example, a simple disc
\begin{center}
$D = \{\varnothing, \{1\}, \{2\}, \{3\}, \{1, 2\}, \{1, 3\},
\{2, 3\}, \{1, 2, 3\}\}$
\end{center}
is a closed simplex, while its boundary
\begin{center}
$\partial D = \{\varnothing, \{1\}, \{2\}, \{3\}, \{1, 2\},
\{1, 3\}, \{2, 3\}\}$
\end{center}
is not, because the triple $\{1,2,3\}$ is not connected.

Examples of simplicial complexes are provided by
triangulations of continuous manifolds.
Note, however, that these should be
{\it triangulations}, not just arbitrary discretizations.
For example, rectangular lattices with hypercubic, rather
than tetrahedra (simplicial) sites are {\it not}
simplicial complexes according to our definition
(because diagonal of a square is not a link).
They can be easily completed to triangulation by adding
diagonals (this, however, can be done in different ways,
since there are two possible choices of diagonal in each
square; independence of the choice of triangulation in this context
is often called {\it flip symmetry}).

We also consider an ordering of vertices in $M$.
Linear ordering is a function $ord(\cdot)$ that labels
vertices $x$ with integer (and even natural) numbers $ord(x)$.
Given any simplex $\sigma$ and vertex $x \in \sigma$, we define
$$ord (x \rightarrow \sigma)$$
to be the number of vertices $y \in \sigma$ with $ord(y) \leq ord(x)$.
For a generic sub-simplex $\sigma' \in \sigma$ 
\be
ord (\sigma' \rightarrow \sigma) =
\sum\limits_{x \in \sigma'} ord (x \rightarrow \sigma)
\ee
Here $|\sigma|$ denotes the dimension of $\sigma$. The number
$\beta_{\sigma' \rightarrow \sigma}
= (-1)^{ord (\sigma' \rightarrow \sigma)} = \pm 1$
defines the parity of embedding of a subset $\sigma' \subset \sigma$ into an ordered set $\sigma$.
It will be used below to specify coefficients in the definitions of various operations.

\subsection{Discrete forms}

Given simplicial complex $M$, a \textit{discrete p-form}
is a function defined on its $p$-simplexes.
$0$-forms are functions on vertices,
while $1$-forms are defined on links, $2$-forms - on triangles, etc. That is a discrete (finite-difference) analogue of smooth forms, because a smooth $p$-form is a function on $p$-dimensional submanifolds. Its value on a particular submanifold is equal to
the integral over this submanifold.

\paragraph{}
Simplexes can be formally added and multiplied by a number. In topology, such linear combinations of simplexes are usually called chains, and linear functions on chains are called co-chains. It is natural to use co-chains to represent discrete forms. However, a discrete form $\omega$ can be treated also a chain:
\begin{center}
$\omega = \sum\limits_{\sigma \in M} \omega (\sigma) \cdot \sigma$
\end{center}
This should not cause any confusion: since a finite-dimensional linear space is isomorphic to its dual space, both points of view are possible. Co-chain description was adopted in s.2, while most of calculations below are done with chains,
but one can always pass to the dual picture.

The wide-spread approach to discrete forms
is the Whitney-form approximation (\cite{prog1}, \cite{prog2}).
In this approach discrete forms are represented by some
collection of smooth forms (Whitney forms)
and discrete operations (derivative, wedge product, etc)
are naturally constructed from their smooth counterparts.
In this paper we study the subject from essentially
different point of view --
we consider discrete forms {\it per se},
as linear combinations of simplexes,
without a reference to any smooth approximation.
Consequently, we use only linear algebra to handle
the discrete forms and operations.

\subsection{Operations}

Given a simplicial complex, denote through $\Omega_* = \oplus_{m=0}^\infty \Omega_m$
the linear space of all discrete forms on it.
$\Omega_*$ is naturally ${\mathbb Z}$-graded into 0-forms $\Omega_0$, 1-forms $\Omega_1$ and so on. The space $\Omega_* \otimes \Omega_*$ is a space of pairs of forms,
$\Omega_* \otimes \Omega_* \otimes \Omega_*$ is a space of triples,
and so on.
The basis in $\Omega_*^{\otimes p}$ is formed by all
possible $p$-uples of simplexes: for example,
$\Omega_* \otimes \Omega_*$ is spanned by pairs of simplexes.

An \textit{operation} with forms is any linear map (tensor)
\begin{center}
$O:\ \ \Omega_*^{\otimes p} \rightarrow \Omega_*^{\otimes q}$
\end{center}
A composition of operations
$$O_{1} \circ O_{2} (\sigma_1 \otimes \sigma_2 \otimes \ldots \otimes \sigma_p) = O_{1} \Big( O_{2} (\sigma_1 \otimes \sigma_2 \otimes \ldots \otimes \sigma_p) \Big)$$
and duality are immediately present. In addition we define a property, inherited from the
structure of simplicial complex: namely, locality.

\subsection{Locality}

Given $p$ simplexes
$\sigma_{1}, \sigma_{2}, \ldots, \sigma_{p}$,
we define their \textit{local envelope}
$$\cup (\sigma_{1}, \sigma_{2},.. \sigma_{p})$$
to be a minimal simplex which contains all of them.
If such a simplex does not exist, the local envelope
is $\varnothing$
(for example, $\cup(\partial D)= \varnothing$).
%\linebreak
For any simplex $\sigma$,
denote through $\Omega_*^{\otimes p} (\sigma)$
a subspace spanned by $p$-uples
with local envelope equal to $\sigma$.
Then the whole space $\Omega_*^{\otimes p}$
is decomposed into subspaces:
\be
\Omega_* ^{\otimes p} =
\bigoplus\limits_{\sigma \in M} \Omega_* ^{\otimes p}(\sigma)
\ee

An operation $O$ is said to be \textit{local}
if it maps $\Omega_*^{\otimes p}$ into $\Omega_*^{\otimes q}$
and never decreases the local envelope:
\be
{\rm local} \ \ \
O:\ \ \Omega_* ^{\otimes p} (\sigma) \mapsto
\bigoplus\limits_{\sigma' \supset \sigma}
\Omega_* ^{\otimes q}(\sigma')
\ee
In a similar way, a \textit{co-local} operation
$O:\ \ \Omega_*^{\otimes p} \rightarrow \Omega_*^{\otimes q}$
is the one that never increases the envelope:
\be
{\rm co-local} \ \ \
O:\ \ \Omega_* ^{\otimes p} (\sigma) \mapsto
\bigoplus\limits_{\sigma' \subset \sigma}
\Omega_* ^{\otimes q}(\sigma')
\ee
Operation that is both local and co-local,
is called \textit{strictly local}.
A strictly local operation
$O:\ \ \Omega_*^{\otimes p} \rightarrow \Omega_*^{\otimes q}$
maps
$\Omega ^{\otimes p} (\sigma)$ into
$\Omega_* ^{\otimes q}(\sigma)$
and hence preserves the local envelope.
Obviously, a composition of (co-)local operations is (co-)local.
Duality converts a local operation into co-local and vice versa.

The meaning of locality becomes clear in the dual picture --
consider discrete forms as functions on simplexes.
Each operation is fully specified, if values of its image on all
simplexes are known.
However, this value on simplex $\sigma$ can depend
on the values of operation's arguments on larger domain of our
simplicial complex which lie beyond $\sigma$.
It is exactly the case when we call the operation non-local.

\paragraph{Localisation.} It is especially important that a non-local operation $O$
still can contain a nonzero local part, because its value $O(e)$ where
$e \in \Omega_* ^{\otimes p} (\sigma)$
can have a nonzero projection onto the local subspace
$\bigoplus\Omega_* ^{\otimes q}(\sigma' \supset \sigma)$.
The same happens to co-local part and to strictly local part.
A generic operation $O$ is a sum of three parts:
$$O = O_{-} + [O] + O_{+}$$

$$[O] = \Pi\big[ \bigoplus\Omega_* ^{\otimes q}(\sigma' = \sigma) \big] \circ O$$

$$O_+ = \Pi\big[ \bigoplus\Omega_* ^{\otimes q}(\sigma' \supset \sigma, \sigma' \neq \sigma) \big] \circ O$$

$$O_- = \Pi\big[ \bigoplus\Omega_* ^{\otimes q}(\sigma' \subset \sigma, \sigma' \neq \sigma) \big] \circ O$$
where $\Pi\big[L\big]$ denotes the projection onto a linear space $L$. From now on $[O]$ denotes the strictly local part of $O$, while $[O] + O_{+}$ and $[O] + O_{-}$
denote the local and co-local parts respectively. Such extraction of a strictly local operation $[O]$ from an originally non-local operation $O$ is called localization of $O$.

\subsection {Discrete $d$ and $\wedge$}

In order to discretise differential forms,
especially important are the analogues of smooth
exterior derivative $d$ and exterior wedge product $\wedge$.
We define them on any ordered simplicial complex.
From now on
\be
\beta_{\sigma' \rightarrow \sigma}
= (-1)^{ord (\sigma' \rightarrow \sigma)}
\ee

\paragraph {Discrete de Rham differential}
$d:\ \ \Omega_* \rightarrow \Omega_*$
is defined as follows:
\be
d (\sigma) = \sum\limits_{x \notin \sigma}
\beta_{x \rightarrow \sigma \cup x} \cdot \cup(\sigma, x)
\ee
Here $\sigma$ is a simplex in $M$, while $x$ - a vertex in $M$. Operator $d$ has a dual (conjugate) operator $\partial$, the co-differential:
\be
\partial (\sigma) = d^{+} (\sigma) =
\sum\limits_{x \in \sigma} \beta_{x \rightarrow \sigma}
\cdot (\sigma \setminus x)
\ee
Discrete Stokes theorem $d \omega (\sigma) = \omega (\partial \sigma)$ is nothing but the duality between $d$ and $\partial$.

In other words, there is a natural pairing between chains and forms:
$$ < \sigma, \omega > = < \omega, \sigma > = \omega (\sigma)$$
Operators $d$ and $\partial$ are conjugate w.r.t this pairing:
$$< \partial(\sigma), \omega > = < \sigma,  d(\omega) >$$

Differential maps a $p$-form into a $(p+1)$-form,
and co-differential converts a $p$-form into a $(p-1)$-form.
From the dual point of view, differential adds one vertex to a simplex
in all possible ways, with appropriate signs,
and co-differential removes one vertex from a simplex
in all possible ways, with appropriate signs.
\paragraph{}
Operations $d$ and $\partial$ can be found in any
combinatorial topology textbook (see, for example, \cite{Pont}),
and have the well-known properties.
First, these operators are nilpotent:
\be
d^2 = \partial ^2 = 0
\ee
This property is most widely used in combinatorial topology:
homology groups of $d$ (``simplicial homologies'')
are applied to the study of topological properties of
simplicial complexes, in particular, of triangulated manifolds.
Second, the anti-commutator
\be
d \partial + \partial d = \bigtriangleup
\ee
is a natural discrete Laplace operator.
Its properties are related to topology of the complex $M$
via Hodge theory.
If $M$ is just a closed simplex,
then $\bigtriangleup$ is very simple:
it is a scalar operator, multiplying any form by
the number of vertices in the simplex, i.e. does not
depend on the choice of the form.
If $M$ is not exactly a simplex, but still has topology of a disc,
then $\bigtriangleup$ is non-degenerate.
Finally, if $M$ has non-trivial topology,
then $\bigtriangleup$ becomes degenerate,
its kernel consists of harmonic $p$-forms
and represents the $p$-th homology of $M$.

\paragraph {Discrete wedge product}

$\wedge: \Omega_* \otimes \Omega_* \rightarrow \Omega_*$
is defined as follows:
\be
\wedge (\sigma, \sigma') =
\dfrac{|\sigma|!|\sigma'|!}{(|\sigma| + |\sigma'| + 1)!}
\,\beta_{\sigma \rightarrow \sigma \cup \sigma'} %\cdot
\beta_{\sigma \cap \sigma' \rightarrow \sigma'} \cdot
\cup(\sigma, \sigma')
\ee
Here $\sigma$ and $\sigma'$ are simplexes in M,
intersecting at a point: $dim (\sigma \cap \sigma') = 0$,
otherwise wedge product is zero.
This product is not as popular in literature
as the discrete de Rham differential, however,
it is indeed the correct discrete realization of exterior product.
First, as a consequence of the intersection rule,
it converts a $p$-form and a $q$-form into a $(p+q)$-form.
Second, as a consequence of the choice of signs $\beta$,
it is skew-symmetric:
\be
\wedge (\omega, \omega') =
(-1)^{|\omega| \cdot |\omega'|} \wedge(\omega', \omega)
\ee
Finally,
it satisfies a graded Leibnitz rule with the differential $d$:
\be
d(\wedge(\omega, \omega')) =
\wedge(d(\omega), \omega') + (-1)^{|\omega|} \cdot
\wedge(\omega, d(\omega'))
\ee
Associativity is the only property of smooth exterior product
that is broken.
However, we do know that instead discrete exterior product
is a part of non-trivial $A_\infty$ structure.
This structure is most natural at the level of simplicial complex,
and associativity arises only in the continuum limit.

%\pagebreak
\section {$A_{\infty}$ structure}

\subsection {The detailed/multilinear form}
The following definition of
$A_{\infty}$ structure or algebra is most widely used
\cite{A8}.

\textit{$A_{\infty}$ structure} on a ${\mathbb Z}$-graded linear space
$\Omega_*$ is a countable collection of multilinear maps
\be
m^{(p\,)}:\ \ \Omega_*^{\times p} \rightarrow \Omega_*,
\ \ \ \ \ p=1,2,\ldots
\ee
of degrees ${\rm deg}\big(m^{(p\,)}\big) = 2 - p\ $
that satisfy the relations
\be
\sum\limits_{k+l = n + 1} \sum\limits_{j=0}\limits^{k-1}
(-1)^{\xi(\!j,l)} m^{(k)}
(\omega_{1},...,\omega_{j},m^{(l)}(\omega_{j+1},\ldots,
\omega_{j+l}),\omega_{j+l+1},\ldots,\omega_{n}) = 0,
\ \ \ \ \ n =1,2,\ldots
\label{A8full}
\ee
The signs factors are defined by
\be
\xi(\!j,l) = (j+1)(l+1) +
l(|\omega_{1}| + ... + |\omega_{j}|)
\label{xisigns}
\ee
where $|\omega|$ is the degree of a form $\omega \in \Omega_*$.
Degree of a $p$-linear map $X$ is defined as
\be
{\rm deg}(X) = |X(\omega_{1}, \ldots , \omega_{p})| -
(|\omega_{1}| + ... + |\omega_{p}|)
\ee

Explicitly, the first three identities read
$$m^{(1)} \Big(m^{(1)} (a)\Big) = 0,$$
$$m^{(1)} \Big( m^{(2)} (a,b)\Big) -
m^{(2)} \Big(m^{(1)} (a),b\Big) -
(-1)^{|a|} m^{(2)} \Big( a, m^{(1)}(b)\Big) = 0,$$
$$m^{(2)} \Big(m^{(2)}(a,b),c\Big) -
m^{(2)} \Big(a,m^{(2)}(b,c)\Big) = $$
$$ = m^{(1)} \Big(m^{(3)}(a,b,c)\Big) +
m^{(3)}\Big(m^{(1)}(a),b,c\Big) +
(-1)^{|a|} m^{(3)}\Big(a,m^{(1)}(b),c\Big) +
(-1)^{|a|} (-1)^{|b|} m^{(3)}\Big(a,b,m^{(1)}(c)\Big)$$

Bilinear map $m^{(2)}: \Omega_* \times \Omega_* \rightarrow \Omega_*$
is usually thought to be a product on $\Omega_*$,
while operator $m^{(1)}: \Omega_* \rightarrow \Omega_*$
is thought to be a differential.
The first $A_{\infty}$-identity implies that $m^{(1)}$
is square-zero, the second one implies that $m^{(1)}$
is a graded derivation of $m^{(2)}$ (Leibnitz rule)
and the third one implies that the product $m^{(2)}$ is
associative up to the higher operation $m^{(3)}$.
Such abstract $A_{\infty}$ structure has a lot of realizations,
one of the most interesting arises in quantum field theory.

However, the physical systems are usually infinite-dimensional. The finite-dimensional $A_{\infty}$ structure on simplicial complexes is more simple
and transparent. In this case, the graded linear space $\Omega_*$ is the space of
discrete forms on a complex,
operator $m^{(1)}$ has a meaning of discrete de Rham derivative
and $m^{(2)}$ is naturally identified with discrete wedge product.

\subsection {The brief/nilpotent form}

In section 2, we introduced a {\it brief} formulation of the
$A_\infty$-structure -- a simplification that happens
if we define $m^{(p)}$ as a linear operation acting not only on
$\Omega_*^{\otimes p}$, but on arbitrary space $\Omega_*^{\otimes q}$
with arbitrary $q$:
$$m^{(p\,)}: \Omega_*^{\otimes q} \rightarrow 0 \ \ \ \ \
{\rm for} \ \ \ q < p$$
$$m^{(p\,)}: \Omega_*^{\otimes q} \rightarrow
\Omega_*^{\otimes q - p + 1} \ \ \ \ \ \ {\rm for} \ \ \ q \geq p$$
The latter action is called a ``lifting rule'' and is defined for $q \geq p$ as follows:
\be
m^{(p\,)} (\omega_{1} \otimes \ldots \otimes \omega_{q}) =
(-1)^{q-1}\sum\limits_{j = 0}\limits^{q-p}
(-1)^{\xi(j,p\,)}
\omega_{1} \otimes \ldots \otimes\omega_{k}
\otimes m^{(p)}
\Big(\omega_{k + 1},\ldots,\omega_{k + p}\Big) \otimes
\omega_{k+p+1} \otimes\ldots
\otimes \omega_{q}
\label{comul}
\ee
where $\xi$-function in the sign factors is given by the same
(\ref{xisigns}).

Actually, in this way we lift our operations to
the total linear Fock-like space
\be
\Omega_*^* = \bigoplus_{p = 0}^{\infty}\Omega_*^{\otimes p}
\ee
(we put $\Omega_*^{\otimes 0} = \mathbb{C}$).
This space is also called a tensor algebra over $\Omega_*$. Note that $\Omega_*^*$ is doubly graded:
there is a conventional degree-grading inherited from $\Omega_*$
plus a tensor-grading: into forms, pairs of forms,
triples of forms and so on. 
By means of lifting,
operations $m^{(p)}$ (previously the poly-linear maps)
act as linear \textit{operators} in this space.
Operators are much easier to handle,
as compared to the multilinear maps.
In these terms the system of constraints (\ref{A8full})
looks substantially better:
$$
\sum_{k + l = n + 1} m^{(k)} \circ m^{(l)}
(\omega_{1} \otimes \ldots \otimes \omega_{n}) = 0
$$
or simply
\be
\sum_{k + l = n + 1} m^{(k)} m^{(l)} = 0
\label{rmk}
\ee
Eq.(\ref{rmk}) is a {\it brief} form of (\ref{A8full}).
\paragraph{}
Explicitly, these constraints read
$$(m^{(1)})^2 = 0$$
\be
m^{(1)} m^{(2)} + m^{(2)} m^{(1)} = 0
\label{explico}
\ee
$$m^{(1)} m^{(3)} + m^{(3)} m^{(1)} + (m^{(2)})^2 = 0$$
and so on.
They can be nicely unified into a single (tensor-)graded constraint:
\be
\big( D \big)^{\circ 2} =
\big( m^{(1)} + m^{(2)} + m^{(3)} + \ldots \big)^{\circ 2} = 0
\ee
In other words, the $A_{\infty}$ structure
on simplicial complex is nothing but a single operator
$D:\ \ \Omega_*^* \rightarrow \Omega_*^*$,
\be
D = m^{(1)} + m^{(2)} + m^{(3)} + \ldots
\ee
satisfying a single nilpotency constraint
\be
D^2 = 0
\label{nilpo}
\ee
The peculiarity is that the linear space $\Omega_*^*$ is graded,
and the single relation (\ref{nilpo}) --
simple and natural at the level of the whole tensor space --
gives rise to sophisticated and numerous constraints
(\ref{A8full}) at the level of graded components.
{\it Brief} formulation removes the grading and
reveals the underlying nilpotent structure.

\subsection{Solving $A_\infty$ relations for $m^{(p)}$}

In what follows
we consider only the special case of $m^{(1)} = d$ and
$m^{(2)} = \wedge$ defined in 4.5.
This restriction is implied by consideration of continuum limit:
for high-resolution triangulations of smooth manifolds,
we want the $A_\infty$ structure to reproduce the smooth algebra
(\ref{nccont}) with increasing precision.
In other words, in this paper we look only at $A_\infty$ structures
with the correct continuum limit.
This natural constraint is very strong:
it fixes $m^{(1)} = d$ and $m^{(2)} = \wedge$ unambiguously.
These and only these two discrete operations reproduce
the smooth operations $d$ and $\wedge$ in the continuum limit,
while all the higher operations $m^{(p)}$ with $p\geq 3$ vanish.

With this choice of $m^{(1)}$ and $m^{(2)}$
the two first $A_\infty$-equations (\ref{explico}) are satisfied:
$$d^2 = 0$$
$$d \wedge + \wedge d = 0$$
The higher $A_\infty$ equations (\ref{rmk}) have
a recursive structure:
\be
d m^{(p\,)} + m^{(p\,)} d = {\cal M}^{(p\,)}
\label{A88}
\ee
where ${\cal M}^{(p\,)}$ is the right hand side
\be
{\cal M}^{(p\,)} =
- \sum_{k = 2}^{p - 1} m^{(k)} m^{(p - k + 1)}
\ee
so that ${\cal M}^{(1)} = 0$, ${\cal M}^{(2)} = 0$,
${\cal M}^{(3)} = -\wedge^2$ etc.
It is natural to solve such equations one by one (recursively),
because ${\cal M}^{(p\,)}$ depends only on the previous operations
$m^{(1)}, \ldots, m^{(p-1)}$.

\subsection{$K$-operator method \label{concon}}

On each step of recursion we need to solve an equation
\be
d m^{(p\,)} + m^{(p\,)} d = {\cal M}^{(p\,)}
\ee
where operation $m^{(p)}$ is the only unknown quantity and
${\cal M}^{(p\,)}$ at the right hand side is already known.
Such equations can be solved by generic ``$K$-operator'' method
\cite{KOP},
which can be roughly described as finding the inverse operator
for differential.
Consider generic equation $d x + x d = y$ and any operator K,
satisfying
\be
d K + K d = 1
\label{K-ident}
\ee
Then both $x = yK$ and $x = Ky$ are solutions,
provided only $dy = yd$ holds,
which is the \textit{consistency condition}.
They are indeed solutions:
$$d(yK) + (yK)d = ydK + yKd = y(dK + Kd) = y$$
$$d(Ky) + (Ky)d = dKy + Kdy = (dK + Kd)y = y$$
We call them the right-hand and left-hand solutions respectively.
In our problem, $x = m^{(p)}$ and $y = {\cal M}^{(p\,)}$.

To begin with, consistency condition
\be
d {\cal M}^{(p\,)} = {\cal M}^{(p\,)} d
\ee
has to be checked.
For $p = 3$ consistency $d \wedge^2 = \wedge^2 d$ follows from
$d \wedge + \wedge d = 0$.
For generic $p \geq 3$ the check is straightforward:
$${\cal M}^{(p\,)} d - d {\cal M}^{(p\,)} =
\sum_{k = 2}^{p - 1} d m^{(k)} m^{(p - k + 1)} - m^{(k)}
m^{(p - k + 1)} d = $$
$$= \sum_{k = 2}^{p - 1} d m^{(k)} m^{(p - k + 1)} +
m^{(k)} d m^{(p - k + 1)} - m^{(k)} {\cal M}^{(p - k + 1)} = $$
$$=\sum_{k = 2}^{p - 1} (d m^{(k)} + m^{(k)} d) m^{(p - k + 1)}
- m^{(k)} {\cal M}^{(p - k + 1)} = %$$ $$ =
\sum_{k = 2}^{p - 1} {\cal M}^{(k)} m^{(p - k + 1)} - m^{(k)}
{\cal M}^{(p - k + 1)} =$$
$$= \sum_{k = 2}^{p - 1} \sum_{q = 2}^{k - 1} m^{(q)}
m^{(k - q + 1)} m^{(p - k + 1)} -\sum_{k = 2}^{p - 1}
\sum_{q = 2}^{p - k} m^{(k)} m^{(q)} m^{(p + 2 - q - k)} = $$
$$= \sum_{k = 2}^{p - 1} \sum_{\alpha = k-1}^{2}
m^{(k - \alpha + 1)} m^{(\alpha)} m^{(p - k + 1)} -
\sum_{\beta = p-1}^{2} \sum_{q = 2}^{\beta - 1}
m^{(p - k + 1)} m^{(q)} m^{(\beta - q + 1)} =$$
$$=\sum_{k = 2}^{p - 1} \sum_{q = 2}^{k-1} \big{(}
m^{(k - q + 1)} m^{(q)} m^{(p - k + 1)} -
m^{(p - k + 1)} m^{(q)} m^{(k - q + 1)} \big{)} = 0$$
The right-hand solution of (\ref{A88}) is
$m^{(p\,)} = {\cal M}^{(p\,)} K$ and the left-hand solution is
$m^{(p\,)} = K {\cal M}^{(p\,)}$.
Note that the freedom is not exhausted by the
single choice between the right and left solutions:
it is possible to independently choose either one or another
at each step of recursion.
If we select only right-hand solutions each time,
a pure right $A_\infty$ structure arises:
$$m^{(1)} = d$$
$$m^{(2)} = \wedge$$
$$m^{(3)} = {\cal M}^{3} K = - \wedge^2 K$$
$$m^{(4)} = {\cal M}^{4} K =
- (\wedge m^{(3)} + m^{(3)} \wedge) K =
(\wedge \wedge^2 K + \wedge^2 K \wedge) K =
\wedge^3 K^2 + \wedge^2 K \wedge K$$
$$m^{(5)} = {\cal M}^{5} K =
- (\wedge m^{(4)} + m^{(4)} \wedge + m^{(3)} m^{(3)}) K =
- \big{(} \wedge ( \wedge^3 K^2 + \wedge^2 K \wedge K) +
(\wedge^3 K^2 + \wedge^2 K \wedge K) \wedge + \wedge^2 K \wedge^2 K
\big{)} K$$
and so on.
In what follows only this pure right-hand construction
will be considered.
A pure left-hand $A_\infty$ structure can be constructed in
a similar way,
however, in section 7 we explain that only the right-hand
construction can produce a local $A_\infty$ structure.
If we are not interested in locality,
then both left or right choices can be made on equal footing.

It is not necessary but often convenient to consider
nilpotent $K$-operators, satisfying
\[
\begin{array}{ccc}
d^2 = 0 \\
\\
d K + K d = 1 \\
\\
K^2 = 0
\end{array}
\]
Such $K$-operators provide especially simple
$A_\infty$ structures:
$$m^{(1)} = d$$
$$m^{(2)} = \wedge$$
$$m^{(3)} = {\cal M}^{3} K = - \wedge^2 K$$
$$m^{(4)} = {\cal M}^{4} K =
- (\wedge m^{(3)} + m^{(3)} \wedge) K = (\wedge \wedge^2 K +
\wedge^2 K \wedge) K = \wedge^2 K \wedge K$$
$$m^{(5)} = {\cal M}^{5} K =
- (\wedge m^{(4)} + m^{(4)} \wedge + m^{(3)} m^{(3)}) K =
- \wedge^2 K \wedge K \wedge K$$
$$\ldots$$
\be
m^{(p)} = (-1)^p \wedge (\wedge K) ^{p - 2}
\label{KNil}
\ee
In this case, the whole series
$D = m^{(1)} + m^{(2)} + m^{(3)} + \ldots$
can be collected into a compact expression.
Indeed,
$$D = d + \wedge - \wedge^2 K + \wedge^2 K \wedge K + \ldots
= d + \wedge \big{(} 1 - \wedge K + \wedge K \wedge K +
\ldots \big{)}$$
Making use of $1/(1 + A) = 1 - A + A^2 - A^3 \ldots$,
it is further simplified:
$$
D = d + \wedge \dfrac{1}{1 + \wedge K} =
\big(d + d \wedge K + \wedge\big)\dfrac{1}{1 + \wedge K} =
\big(d + d \wedge K + \wedge (d K + K d)\big)
\dfrac{1}{1 + \wedge K} =
$$
\be
= \big(d + \wedge K d\big)\dfrac{1}{1 + \wedge K} =
(1 + \wedge K)\ d\ (1 + \wedge K)^{-1}
\label{Conj}
\ee
It turns out that the $A_\infty$ nilpotent operator
$D$ is a conjugation of the \textit{bare} nilpotent operator $d$:
$D = U d U^{-1}$.
The transformation map $U$ is equal to
\be
U = 1 + \wedge K
\ee
in this particular case (it would be $U = 1 + K \wedge$
if we choose   the pure left-hand solutions).
Of course, conjugation preserves nilpotency:
\be
D^2 = U d U^{-1} U d U^{-1} = U d^2 U^{-1} = 0
\ee

\section{The first (naive) solution}

There is a natural way to construct a $K$-operator from
\textit{any} operator $\alpha:\ \Omega_*^* \rightarrow \Omega_*^*$.
Introduce
$$
d \alpha + \alpha d = \bigtriangleup_{\alpha}
$$
Nilpotency $d^2=0$ of $d$ implies that
$$
d\bigtriangleup_{\alpha} = \bigtriangleup_{\alpha}d
$$
If this $\bigtriangleup_{\alpha}$ is invertible,
then the operator
$$
K_{\alpha} = \alpha \cdot \bigtriangleup^{-1}_{\alpha}
$$
satisfies the required identity
\be
d K_{\alpha} + K_{\alpha} d = 1
\ee
Indeed,
$$d K_{\alpha} + K_{\alpha} d =
d \alpha \bigtriangleup^{-1}_{\alpha}
+ \alpha  \bigtriangleup^{-1}_{\alpha} d =
d \alpha \bigtriangleup^{-1}_{\alpha} +
\alpha d \bigtriangleup^{-1}_{\alpha} =
(d \alpha + \alpha d) \bigtriangleup^{-1}_{\alpha} =
\bigtriangleup_{\alpha} \bigtriangleup^{-1}_{\alpha} = 1$$
The same is true also for $\tilde K_\alpha =
\bigtriangleup_{\alpha}^{-1}\alpha$.
Of course, $\bigtriangleup_{\alpha}$ is rarely invertible.
For arbitrary $\alpha$ this operator does not
respect grading and $\bigtriangleup_{\alpha}$ acts from
a sub-space of a given grading into another.
Even if the sub-space remains intact, there can be zero modes.

The simplest nontrivial operator $\alpha$ one can choose is
$\alpha = \partial$, so that
$\Delta_\partial = d \partial + \partial d\ = \Delta$
is the ordinary Laplace operator, which respects all gradings.
If $M$ is just a closed simplex, then $\bigtriangleup$
is most simple: it is a scalar operator,
multiplying any form by the number of vertices in the simplex.
Therefore $\alpha = \partial$ is especially convenient for
calculations on a simplex.
However, it is also of use on generic simplicial complexes.

If we choose $\alpha = \partial$, then, since $\partial^2=0$,
$K = K_{\partial} = \partial \cdot \bigtriangleup^{-1} =
\bigtriangleup^{-1}\cdot \partial $
and occasionally
$$K^2 = \partial \bigtriangleup^{-1} \partial \bigtriangleup^{-1}
= \partial \partial \bigtriangleup^{-1} \bigtriangleup^{-1} = 0$$
Therefore, according to (\ref{KNil}), we obtain the
following (left-hand) $A_{\infty}$ structure:
\be
\boxed{
\left\{
\begin{array}{ll}
m^{(1)} = d,\\
m^{(p)} = (-1)^p \Big( \partial\
({d \partial + \partial d})^{-1} \wedge \Big) ^{p - 2} \wedge & p \geq 2\\
\end{array}
\right.}
\label{SOL1}
\ee
Note that this solution is as explicit as a formula can be!
The differentials $d$, $\partial$ and the wedge multiplication
$\wedge$ are defined in section 4.5 in pure combinatorial
terms, and such operations $m^{(p)}$ can be calculated easily
on various simplicial complexes.
\pagebreak
\section{The second (local) solution}

The solution (\ref{SOL1}) is just an illustration
of how the method works:
we pick a $K$-operator and immediately obtain a solution.
Another $K$-operator would produce another $A_{\infty}$ structure.
Actually, the $K$-operator method is capable of solving more
complicated problems.
It can be used to construct not just any but specific
$A_\infty$ structures, with given properties, for example,
{\it locality}.
In this section we are going to construct a local
$A_\infty$ structure.
An $A_\infty$ structure is said to be local,
if all its operations
$m^{(p)}:\ \Omega_*^{\otimes p} \rightarrow \Omega_*$ are local.
This is a strong additional requirement.

Local $A_\infty$ structure can be built with the help of a
\textit{local K-operator}.
Note that neither $\partial$ nor $K_{\partial}$ are local,
and the solution (\ref{SOL1}) is not local as well.
To construct a local solution,
let us consider some important properties of the discrete differential.

\subsection {Localisation of $d$ and $\partial$}

\paragraph{$\bullet$} In section 4.4 we have outlined that any operation can be decomposed into three operations with definite locality properties. This decomposition can be applied to operations $d$ and $\partial$:

$$d = d_- + d_+ + [d]$$

$$\partial = \partial_- + \partial_+ + [\partial]$$

\paragraph{$\bullet$} We are especially interested in describing the strictly local parts [$d$] and [$\partial$]. First of all, on the space of forms $\Omega_*$ the local envelope is trivial: for any simplex $\sigma$, $\cup (\sigma) = \sigma$.
Differential $d(\sigma)$ increases $\sigma$ and the local
envelope by adding a vertex from beyond $\sigma$.
Consequently, the strictly local part of $d$ vanishes.
The same is true for $\partial$:
it decreases the local envelope by removing a vertex. Therefore we obtain
$$[d]\Big{\vert}_{\Omega_*} = 0$$
$$[\partial] \Big{\vert}_{\Omega_*} = 0$$
Remarkably, operations [$d$] and [$\partial$] become non-vanishing
on extended spaces $\Omega_* ^{\otimes p}$  with $p\geq 2$.
Operations $d$ and $\partial$ are defined on these spaces
via lifting (\ref{comul}).
For example, on the space of pairs $\Omega_* \otimes \Omega_*$
the differential act on the first simplex,
and then on the second one:
$$d (a \otimes b) = - d(a) \otimes b - (-1)^{|a|} a \otimes d(b)$$
$$\partial (a \otimes b) =
- \partial(a) \otimes b - (-1)^{|a|} a \otimes \partial(b)$$
There can be vertices beyond $a$ that belong to $b$
-- adding such a vertex to $a$ does not change the envelope.
This is what allows $d(a \otimes b)$ to contain a strictly local
(envelope-preserving) part.
There can be vertices that belong to both $a$ and $b$,
removing such a vertex from $a$ does not change the envelope.
Therefore $\partial(a \otimes b)$ contains a strictly local
(envelope-preserving) part. \footnote{In particular, if simplices $\sigma_1, \sigma_2, \ldots, \sigma_p$ have no pairwise intersection points, then $[\partial](\sigma_1 \otimes \sigma_2 \otimes \ldots \otimes \sigma_p) = 0$}

Consider the following minimal example:
the complex is a $1$-simplex
$M = \{\varnothing, \{1\}, \{2\}, \{1,2\}\}$, then
$$\partial \Big(\{1\} \otimes \{1,2\}\Big) =
- \varnothing \otimes \{1,2\}
- \{1\} \otimes \{2\} + \{1\} \otimes \{1\}$$
and
$$[\partial] \Big(\{1\} \otimes \{1,2\}\Big) =
- \varnothing \otimes \{1,2\} - \{1\} \otimes \{2\}$$
The term $\{1\} \otimes \{1\}$ is non-local in this case,
because its local envelope $\{1\}$ does not contain $\{1,2\}$,
which is the local envelope of the argument.
One can see, that localization is nothing but omitting
non-local terms.

\paragraph{$\bullet$} Differential $d$ is a local operation on the space $\Omega_*$ of single forms. However, on the whole Fock space $\Omega^*_*$ it contains a non-local part. The non-local part of $d$ (operation $d_-$) is not identically zero for a very special reason: because the value $d(\sigma_1, \sigma_2, \ldots, \sigma_p)$ can have empty envelope, even if the argument $\sigma_1, \sigma_2, \ldots, \sigma_p$ had envelope $\sigma \neq \varnothing$. For example, consider a pair $\{1,2\} \otimes \{2\}$ on the simplicial complex $M = \{\varnothing, \{1\}, \{2\}, \{3\}, \{1,2\}, \{2,3\}\}$. This pair has a non-empty envelope $\{1,2\}$. At the same time
$$d (\{1,2\} \otimes \{2\}) = \{1,2\} \otimes \{2,3\} - \{1,2\} \otimes \{1,2\}$$
One can see that the right-hand side contains a nonlocal term: $\cup \big(\{1,2\} \otimes \{2,3\}\big) = \varnothing$.

\paragraph{$\bullet$} In other words, the image of operation $d_-$ always has empty envelope:
$$d_-: \Omega^p_* \rightarrow \Omega^p_*(\varnothing)$$
Due to this fact, the following equality holds: $[d^2] = [d]^2 = 0$. Indeed,
$$[d d] = [ ([d] + d_- + d_+)\cdot([d] + d_- + d_+)] =  [ [d] [d] + d_- d_+ + d_+ d_- ] = [d]^2 + [d_- d_+ + d_+ d_-]$$
Using $\cup (d_- ( \omega ) ) = 0 \ \ \forall \omega$, it is possible to check that the last term $[d_- d_+ + d_+ d_-]$ is zero.

\paragraph{$\bullet$} Finally, let us prove an important identity
\be
d_+ [\partial] + [\partial] d_+ = 0 \label{plusminus}
\ee
We are going to use this identity later. The proof is quite trivial, if the complex $M$ is a simplex. In this case, according to (\ref{Lasi}), the discrete Laplase operator $d \partial + \partial d$ is a scalar operator (proportional to unity) on the space of forms. According to the lifting rule (\ref{comuLa}), it is also scalar on the whole Fock space. It is clear that $ \Delta_+ $ should vanish then, because any scalar operator is strictly local.
$$\Delta_+ = (d \partial + \partial d)+ = d_+ [\partial] + [\partial] d_+ = 0$$
If the complex $M$ is not a simplex exactly, select any $p$-uple of simplexes $e$ with local envelope equal to $\sigma$. Operator $d_+ [\partial] + [\partial] d_+$, should increase the envelope of $e$ by one vertex exactly. Let us consider all such simplexes of $M$ (simplexes, that contain $\sigma$ and one more vertex). The restriction of $d_+ [\partial] + [\partial] d_+$ on any of those simplexes is zero (as stated in the first part of a proof). Therefore the value $(d_+ [\partial] + [\partial] d_+) (e)$ is zero.

\subsection{Operator $\bigtriangleup_{loc} = [d] [\partial] + [\partial] [d]$}

In what follows the central role will be played
by local Laplace operator
\be
\bigtriangleup_{loc} = [d] [\partial] + [\partial] [d]
\ee
On any simplicial complex, irrespective of its topology,
this operator has interesting spectral properties. Consider the Fock space $\Omega^*_* = \oplus\ \Omega_* ^{\otimes p}$. The basis in this space is formed by all finite collections of simplexes
$$
e = \Big(\sigma_{1}, \ldots, \sigma_{p}\Big)
$$
For each basis vector $e$, define $n(e)$ to be the number of
vertices in $\cup (e)$, and $k(e)$ to be the number of
\textit{free} vertices in $\cup (e)$.
A vertex is said to be free, if it belongs to one and only
one simplex from $e$. For example, if $e = \Big(\{1\},\{12\},\{13\}\Big)$, then
$n(e) = 3$, $k(e) = 2$ and $\tau = \{2,3\}$ is the set of free points.

Denote through $\Omega_*^{\otimes p}(\sigma , \tau)$ the subspace,
spanned by collections of $p$ simplexes with given local envelope
$\sigma$ and set of free points $\tau \subset \sigma$.
It is convenient to study the properties of $\bigtriangleup_{loc}$
on these subspaces separately, because the operator does not mix them. The properties of $\bigtriangleup_{loc}$ depend on the numbers $n = |\sigma|$ and $k = |\tau|$. Let us describe them in detail.

\subsection{The eigenvalues}
On the $p^k$-dimensional linear space
$\Omega_*^{\otimes p}(\sigma , \tau)$
where $n = |\sigma|$ and $k = |\tau|$ operator
[$d$][$\partial$] + [$\partial$][$d$]
is self-adjoint, so in the basis of eigenvectors it is diagonal
(no Jordan cells can arise).
The eigenvectors can be written as
$ \mid n_1, n_2, \ldots, n_k\rangle $
with $0 \leq n_i \leq (p - 1)$, and the eigenvalues are
$$\bigtriangleup_{loc} \mid n_1, n_2, \ldots, n_k \rangle\ =
p(n - k + \theta(n_1) + \theta(n_2) + \ldots + \theta(n_k))
\mid n_1, n_2, \ldots, n_k \rangle$$
Here $\theta(x) = \left\{
\begin{array}{c}
1, x > 0  \\
0, x \leq 0 \\
\end{array}
\right.$
It follows that there are $k + 1$ different eigenvalues
\be
\lambda_i = (n - i)p, \ \ \ \ i = 0, \ldots, k
\label{Eigenval}
\ee
with corresponding multiplicities
\be
N_i = C^{i}_{k} (p-1)^i, \ \ \ \ i = 0, \ldots, k
\ee
Obviously, $\sum_{i=0}^k N_i = p^k$. The multiplicity $N_i$ denotes the dimension of i-th eigenspace, that corresponds to the eigenvalue $\lambda_i$. These multiplicities can be written as tables of ``generalized binomial coefficients''
$$p = 2: \left\{
\begin{array}{llllllll}
k = 0: \ \ 1 = 1 & & & &\\
k = 1: \ \ 2 = 1 + 1 & & &\\
k = 2: \ \ 4 = 1 + 2 + 1 & &\\
k = 3: \ \ 8 = 1 + 3 + 3 + 1 &\\
k = 4: \ \ 16 = 1 + 4 + 6 + 4 + 1\\
\end{array} \right.$$
$$p = 3: \left\{
\begin{array}{llllllll}
k = 0: \ \ 1 = 1 & & & &\\
k = 1: \ \ 3 = 1 + 2 & & &\\
k = 2: \ \ 9 = 1 + 4 + 4 & &\\
k = 3: \ \ 27 = 1 + 6 + 12 + 8 &\\
k = 4: \ \ 81 = 1 + 8 + 24 + 32 + 16\\
\end{array}
\right.$$
$$p = 4: \left\{
\begin{array}{llllllll}
k = 0: \ \ 1 = 1 & & & &\\
k = 1: \ \ 4 = 1 + 3 & & &\\
k = 2: \ \ 16 = 1 + 6 + 9 & &\\
k = 3: \ \ 64 = 1 + 9 + 27 + 27 &\\
k = 4: \ \ 256 = 1 + 12 + 54 + 108 + 81\\
\end{array}
\right.
$$
and so on. Note that multiplicities (decomposition of $\Omega_*^{\otimes p}(\sigma , \tau)$ into eigenspaces) depend on the number of free points only.

\subsection{The eigenvectors}
The eigenvectors are given by the following explicit construction.
Suppose first that all the points are free, i.e. $n = k$.
In the simplest case $p = p$, $k = 1$, $n = 1$
the space consists of $p$-uples with only one point
which is itself a free point.
Then in the natural basis
$$(1,0,0,\ldots,0) =
\{pt\} \otimes \varnothing \otimes \ldots \otimes \varnothing$$
$$(0,1,0,\ldots,0) =
\varnothing \otimes \{pt\} \otimes \ldots \otimes \varnothing$$
$$\ldots$$
$$(0,0,0,\ldots,1) =
\varnothing \otimes \varnothing \otimes \ldots
\otimes \{pt\} $$
it is straightforward to prove that the
local Laplace operator is
represented by
the $p \times p$ matrix
\be
\bigtriangleup_{loc} = [d][\partial] + [\partial][d] =
\left(\begin{array}{ccccccccccccc}
p-1 & +1 & -1 & +1 \ldots \\
+1 & p-1 & +1 & -1 \ldots \\
-1 & +1 & p-1 & +1 \ldots \\
\ldots & \ldots & \ldots & \ldots &
\end{array}
\right) =
\label{lapla}
\ee

\be
= p \cdot I -
\left(
\begin{array}{ccccccccccccc}
+1 & 0 & 0 & 0 \ldots \\
0 & -1 & 0 & 0 \ldots \\
0 &  0  & +1 & 0 \ldots \\
\ldots & \ldots & \ldots & \ldots &
\end{array}
\right)
\cdot
\left(
\begin{array}{ccccccccccccc}
+1 & +1 & +1 &+1 \ldots \\
+1 & +1 & +1 &+1 \ldots \\
+1 & +1 & +1 &+1 \ldots \\
\ldots & \ldots & \ldots & \ldots &
\end{array}
\right)
\cdot
\left(
\begin{array}{ccccccccccccc}
+1 & 0 & 0 & 0 \ldots \\
0 & -1 & 0 & 0 \ldots \\
0 &  0  & +1 & 0 \ldots \\
\ldots & \ldots & \ldots & \ldots &
\end{array}
\right) ^{-1}
\ee
where $I$ is the unity matrix.
The matrix (\ref{lapla}) has two eigenvalues, $0$ and $p$,
with the multiplicities $1$ and $p-1$ respectively. Denote
\be
E = \left(
\begin{array}{ccccccccccccc}
+1 & 0 & 0 & 0 \ldots \\
0 & -1 & 0 & 0 \ldots \\
0 &  0  & +1 & 0 \ldots \\
\ldots & \ldots & \ldots & \ldots &
\end{array}
\right)
\ee
The eigenvector $|0\rangle = E \cdot (1,1,1,1\ldots)$
corresponds to $\lambda = 0$,
while $|1\rangle, \ldots, |p-1\rangle$
-- a possible choice is
$$|1\rangle = E \cdot (1,-1,0,0,\ldots)$$
$$|2\rangle = E \cdot (0,1,-1,0,\ldots)$$
$$\ldots$$ correspond to $\lambda = p$.

Using the eigenvectors with one free point, it is possible to construct the generic eigenvectors, for $k \geq 1$ and with some non-free points $e = \sigma_{1} \otimes \ldots \otimes \sigma_{p}$, $k(e) = 0$. The explicit construction is
\be
\mid n_1, n_2, \ldots, n_k\rangle =
e \star |n_1\rangle \star |n_2\rangle \star |n_3\rangle
\star \ldots \star |n_k\rangle
\label{Eigenv}
\ee
The $\star$ operation here attaches one collection of simplexes to another:
$$e_1 \star e_2 = \big(\sigma_1 \otimes \sigma_2 \otimes
\ldots \otimes \sigma_p\big) \star
\big(\sigma'_1 \otimes \sigma'_2 \otimes \ldots
\otimes \sigma'_p\big) =
(-1)^\chi \cup(\sigma_1, \sigma'_1)
\otimes \cup(\sigma_2, \sigma'_2) \otimes \ldots
\otimes \cup(\sigma_p, \sigma'_p)$$
so that
$$\bigtriangleup_{loc} (e_1 \star e_2) =
(\bigtriangleup_{loc} e_1) \star e_2 + e_1 \star (\bigtriangleup_{loc} e_2) \ \ \ \ \forall e_1,e_2$$
This construction explains all the spectral properties. In every particular case, all of them can be easily checked by using only combinatorial definition of $d$, $\partial$ and localization [$\cdot$]. However, some technical details of a generic proof (mainly concerning the sign factors $\chi$) are still missing. This sign $\chi$ depends on various factors including the number and configuration of non-free points. We have not found yet any convenient formula to express it.

\subsection{A local $K$-operator}

Given a strictly local operation $[\partial]$, it is possible to construct a $K$-operator $$K_{[\partial]} = \dfrac{[\partial]}{d [\partial] + [\partial] d}$$
following the procedure from section 6. This operator satisfies the main $K$-operator identity (\ref{K-ident}), however, it is still not local, because $d$ itself has a non-local part (see s.7.1). To solve the problem, we introduce
\be
[K] = \dfrac{[\partial]}{[d] [\partial] + [\partial] [d]}
\label{locaK}
\ee
instead of $K_{[\partial]}$. \footnote{ Zero modes of $[d] [\partial] + [\partial] [d]$ could spoil such a construction. According to (\ref{Eigenval}), zero eigenvalues exist iff $n = k$ -- when all $n$ points are free. In such subspaces the inverse operator $([d][\partial] + [\partial][d])^{-1}$ is not defined unambiguously. However, this unambiguity affects nothing, because $[\partial]$ is zero on such spaces. } This operator is obviously strictly local, but it does not satisfy (\ref{K-ident}). It satisfies only
$$[d] [K] + [K] [d] = 1$$
or equvalently
$$d [K] + [K] d = ([d] + d_+ + d_-) [K] + [K] ([d] + d_+ + d_-) = 1 + (d_+ [K] + [K] d_+) + (d_- [K] + [K] d_-)$$
The second term $d_+ [K] + [K] d_+$ is equal to zero, according to (\ref{plusminus}) and the definition of $[K]$. Unfortunately, the last term $\epsilon = d_- [K] + [K] d_-$ is not exactly zero.
$$d [K] + [K] d = 1 + \epsilon$$
Nevertheless, operator $[K]$ can be used to solve $A_{\infty}$ equations. To prove that ${\cal M}^{(p\,)} [K]$ is a solution, we use the property $\cup(d_-(e)) = \varnothing$ (see s.7.1.). This property implies that composition of $\epsilon$ with any strictly local operation vanishes. In particular, all operations ${\cal M}^{(p\,)}$ are strictly local on $\Omega_*^{\otimes p}$.

The following cancellation of non-locality happens:
$${\cal M}^{(p\,)} \circ \epsilon \big( \sigma_1, \sigma_2, \ldots, \sigma_p \big) = 0$$
Hence $\epsilon$ does not affect the results:
$${\cal M}^{(p\,)} \circ (d [K] + [K] d) = {\cal M}^{(p\,)} + {\cal M}^{(p\,)} \circ \epsilon = {\cal M}^{(p\,)}$$
One can see that ${\cal M}^{(p\,)} [K]$ is indeed a solution:
$$d {\cal M}^{(p\,)} [K] + {\cal M}^{(p\,)} [K] d = {\cal M}^{(p\,)}$$

This cancelation of non-local contributions is one of the main reasons
why the right-hand $K$-operator construction (see 5.4) is used.
Another reason is that operators $[\partial]$ and $[K]$ vanish on the space
$\Omega_*$ of single forms,
so that in the left-hand construction they would be useless:
$$[K] \circ {\cal M}^{(p\,)} \big( \sigma_1, \sigma_2, \ldots, \sigma_p \big) = 0$$

\subsection{Inversion of $\bigtriangleup_{loc}$} Operator $[K]$ includes the inverse of local Laplace operator $\bigtriangleup_{loc}$. We know everything about eigenvectors and eigenvalues of $\bigtriangleup_{loc}$, therefore it is easily invertible: on each subspace $\Omega_*^{\otimes p}(\sigma , \tau)$, where $n = |\sigma|$ and $k = |\tau|$ the identity
\be
(\bigtriangleup_{loc} - \lambda_0 I) \cdot (\bigtriangleup_{loc} -
\lambda_1 I)\cdot \ldots \cdot(\bigtriangleup_{loc} - \lambda_k I)
= 0
\ee
holds, because it is the spectral equation for a
finite-size matrix $\bigtriangleup_{loc}$.
It follows that
\be
\bigtriangleup_{loc}^{-1} = \alpha_0 + \alpha_1 \bigtriangleup_{loc}
+ \ldots + \alpha_k \bigtriangleup_{loc}^k
\label{Spectral}
\ee
Coefficients $\alpha$ are expressed through the eigenvalues
$\lambda$ by Vieta formula:
$$\alpha_{k-i} = (-1)^{i} \cdot S_i (\lambda_0, \lambda_1 \ldots \lambda_k) / S_{k+1} (\lambda_0, \lambda_1 \ldots \lambda_k) $$
where $S_i$ denotes the i-th elementary symmetric polinomial.

This makes the calculation of inverse Laplace operator
straightforward.
However, evaluation of subsequent powers
$\bigtriangleup_{loc}^k$ becomes more and more complicated
when the number of free points is large.
In order to simplify calculations, one can make use of
the eigenvectors and multiplicities and derive an explicit
formula for the action of $\bigtriangleup_{loc}^{-1}$ on any
$p$-uple of simplexes $e$, in terms of $n(e), k(e)$ and $p$.
We present this formula in the Appendix.

\subsection{A local $A_{\infty}$ structure  \label{dco}}

The $K$-operator $[K]$ is nilpotent: it satisfies $[K]^2 = 0$.
Therefore, according to (\ref{KNil}) we obtain an $A_{\infty}$
structure
\be
\boxed{
\left\{
\begin{array}{ll}
m^{(1)} = d,\\
m^{(p\,)} = (-1)^p \wedge \Big( \wedge [\partial]
\ ([d][\partial] + [\partial] [d])^{-1}\Big)^{p - 2}
& p \geq 2\\
\end{array}
\right.
}
\label{ZOL2}
\ee
This $A_{\infty}$ structure is local (it is made of local
operations only) and has a correct continuum limit
(due to correct choice of the first discrete operations $d$
and $\wedge$).
Like (\ref{SOL1}) eq.(\ref{ZOL2}) is an absolutely explicit formula.
In the Appendix below we calculate
a few first operations $m^{(p\,)}$ with the help of (\ref{ZOL2}).
In the simplest cases this calculation can be even done by hand.
In generic case, starting actually from $p = 4$,
this becomes a complicated calculation
and computer facilities can be used.
However, this matters only if one is interested
in explicit expressions for operations.

From pure theoretical point of view most important is that
we obtained a nilpotent $A_{\infty}$ operator
\be
D = (1 + \wedge [K])\ d\ (1 + \wedge [K])^{-1}
\ee

\begin{center}

\end{center}

As already mentioned in s.5.2, this operator is not a
homogenious element w.r.t. the tensor-grading.
It can be decomposed into the graded components:
\be
D = d + \wedge + m^{(3)} + \ldots
\ee
These graded components represent particular
$A_{\infty}$ operations.

\section{Conclusion}

The local $A_\infty$ structure (\ref{ZOL2})
is a particular example
of a wide class of local $A_\infty$ structures.
In other words, equation (\ref{A8''}) has (infinitely) many
solutions and even more, infinitely many local solutions.
A natural claim is that the $A_\infty$ structures
are in one-to-one correspondence with
nilpotent operators on the Fock space $\Omega^*_*$.

Higher nilpotent analogues of $A_\infty$ structure should exist,
being the $\wedge_n$-induced deformations of the $d_n$ operators.
The relevant operators $D_n$ on the Fock space $\Omega^*_*$ satisfy
$$D_n^n = 0$$
For explicit description of the local $A_\infty^{(n)}$
structures see \cite{hini}.

\pagebreak
\section{Appendix}

In this section we present a number of explicit calculations with discrete forms. They should be regarded as examples and illustrations.

\subsection{Simplicial complexes}

\paragraph{One-dimensional complexes:}
\paragraph{}
$M = \{\varnothing,\{1\}, \{2\}, \{1,2\}\}$ - a closed link (1-simplex) \ \ \begin{picture}(60,10)
\put(10,5){\line(1,0){40}}
\put(10,5){\circle*{5}}
\put(50,5){\circle*{5}}
\put(10,10){1} \put(50,10){2} \put(25,-5){12}
\end{picture}
\paragraph{}
$M = \{\varnothing,\{1\}, \{2\}, \{3\}, \{4\}, \{1,2\}, \{2,3\}, \{2,4\}\}$ \ \
\begin{picture}(100,50)
\put(10,-15){\line(1,0){40}}
\put(50,-15){\line(1,0){40}}
\put(50,-15){\line(0,1){40}}
\put(10,-15){\circle*{6}}
\put(50,-15){\circle*{6}}
\put(90,-15){\circle*{6}}
\put(50,25){\circle*{6}}
\put(7,-10){1} \put(53,-10){2} \put(90,-10){4} \put(53,29){3}
\put(25,-25){12} \put(65,-25){24} \put(53,12){23}
\end{picture}

\paragraph{}
$M = \{\varnothing,\{1\}, \{2\}, \{3\}, \{4\}, \{1,2\}, \{1,3\}, \{2,3\}, \{2,4\}\}$ \ \ \begin{picture}(100,50)
\put(10,-15){\line(1,0){40}}
\put(50,-15){\line(1,0){40}}
\put(50,-15){\line(0,1){40}}
\put(10,-15){\line(1,1){40}}
\put(10,-15){\circle*{6}}
\put(50,-15){\circle*{6}}
\put(90,-15){\circle*{6}}
\put(50,25){\circle*{6}}
\put(7,-10){1} \put(53,-10){2} \put(90,-10){4} \put(53,29){3}
\put(25,-25){12} \put(65,-25){24} \put(53,12){23}
\put(20,12){13}
\end{picture}

\bigskip
\paragraph{Two-dimensional complexes:}
\paragraph{}

$M = \{\varnothing,\{1\}, \{2\}, \{3\}, \{1,2\}, \{1,3\}, \{2,3\}, \{1,2,3\}\}$ - a closed triangle (2-simplex)
\paragraph{}
$M = \{\varnothing,\{1\}, \{2\}, \{3\}, \{4\}, \{1,2\}, \{1,3\}, \{1,4\},\{2,3\}, \{2,4\}, \{3,4\}, \{1,2,3\}, \{1,2,4\}, \{1,3,4\}, \{2,3,4\}\}$

- a 2-sphere (a boundary of 3-simplex)

\subsection{Discrete forms and operations: chain description}

On the space of discrete forms various operations are defined. The most important are: the discrete differential $d$ and discrete wedge product $\wedge$, defined in s.4.5. Let us present some explicit calculations with these operations in the chain picture.

\paragraph{The differential.} Consider a complex $M_1 = \{\varnothing, \{1\}, \{2\}, \{3\}, \{1,2\}, \{1,3\}, \{2,3\}, \{1,2,3\}\}$ - a simple 2-disc. A discrete form $f \in \Omega_*(M)$ is a linear combination of these simplexes:

$$f = f_{\varnothing} \cdot \varnothing + f_{1} \cdot \{1\} + f_{2} \cdot \{2\} + \ldots + f_{123} \cdot \{1,2,3\}$$

Differential $d$ acts on $\Omega_*(M)$ as follows:
\begin{center}
$d (\varnothing) = \{1\} + \{2\} + \{3\}$
\end{center}

\begin{center}
$d (\{1\}) = - \{1,2\} - \{1,3\}$

$d (\{2\}) = \{1,2\} - \{2,3\}$

$d (\{3\}) = \{1,3\} + \{2,3\}$
\end{center}

\begin{center}
$d (\{1,2\}) = \{1,2,3\}$

$d (\{1,3\}) = - \{1,2,3\}$

$d (\{2,3\}) = \{1,2,3\}$
\end{center}

\begin{center}
$d (\{1,2,3\}) = 0$
\end{center}

It is easy to check that $d^2 = 0$ is true.
\paragraph{}
The dual operator $\partial$ acts on $\Omega_*(M)$ as follows:

\begin{center}
$\partial (\varnothing) = 0$
\end{center}

\begin{center}
$\partial (\{1\}) = \partial (\{2\}) = \partial (\{3\}) = \varnothing$
\end{center}

\begin{center}
$\partial (\{1,2\}) = \{2\} - \{1\}$

$\partial (\{1,3\}) = \{3\} - \{1\}$

$\partial (\{2,3\}) = \{3\} - \{2\}$
\end{center}

\begin{center}
$\partial (\{1,2,3\}) = \{1,2\} - \{1,3\} + \{2,3\}$
\end{center}

It may be convenient to represent the linear operators with matrices. Obviously, the matrix of $\partial$ is equal to the transposed matrix of $d$.

\begin{center}
on 0-forms \ \ \ \ \ $d = \left[
\begin{array}{cccccc}
-1 & 1 & 0\\
-1 & 0 & 1 \\
0 & -1 & 1 \\
\end{array}
\right]
\ \ \ \ \ \partial = \left[
\begin{array}{cccccc}
1 & 1 & 1 \\
\end{array}
\right]$
\end{center}

\begin{center}
on 1-forms \ \ \ \ $d = \left[
\begin{array}{cccccc}
1 & -1 & 1\\
\end{array}
\right]
\ \ \ \ \ \partial = \left[
\begin{array}{cccccc}
-1 & -1 & 0\\
1 & 0 & -1 \\
0 & 1 & 1 \\
\end{array}
\right]$
\end{center}

\begin{center}
on 2-forms \ \ $d = \left[
\begin{array}{cccccc}
0\\
\end{array}
\right]
\ \ \ \partial = \left[
\begin{array}{cccccc}
1\\
-1\\
1\\
\end{array}
\right]$
\end{center}
It is easy to check that the Laplace operator $\bigtriangleup = d \partial + \partial d = 3 \cdot id$.

\paragraph{}
Consider the second example: $M_2 = \{\varnothing,\{1\}, \{2\}, \{3\}, \{4\}, \{1,2\}, \{2,3\}, \{2,4\}\}$. Here we have

\begin{center}
$d (\varnothing) = \{1\} + \{2\} + \{3\} + \{4\}$
\end{center}

\begin{center}
$d (\{1\}) = - \{1,2\}$

$d (\{2\}) = \{1,2\} - \{2,3\} - \{2,4\} $

$d (\{3\}) = \{2,3\}$

$d (\{4\}) = \{2,4\}$
\end{center}

\begin{center}
$d (\{1,2\}) = d (\{2,3\}) = d (\{2,4\}) = 0$
\end{center}
Pay attention to the changes as compared to $M_1$. They are due to topology of $M_2$: the point $\{2\}$ has 3 incoming links, while other points have only one.

\begin{center}
$\partial (\varnothing) = 0$
\end{center}

\begin{center}
$\partial (\{1\}) =  \partial (\{2\}) = \partial (\{3\}) = \partial (\{4\}) = \varnothing$
\end{center}

\begin{center}
$\partial (\{1,2\}) = \{2\} - \{1\}$

$\partial (\{2,3\}) = \{3\} - \{2\}$

$\partial (\{2,4\}) = \{4\} - \{2\}$
\end{center}

Laplace operator becomes non-trivial,

\begin{center}
On 0-forms \ \ $\bigtriangleup = \left[
\begin{array}{cccccc}
2 & 0 & 1 & 1\\
0 & 4 & 0 & 0\\
1 & 0 & 2 & 1\\
1 & 0 & 1 & 2\\
\end{array}
\right]$
\end{center}

\begin{center}
On 1-forms \ \ $\bigtriangleup = \left[
\begin{array}{cccccc}
2 & -1 & -1 \\
-1 & 2 & 1 \\
-1 & 1 & 2 \\
\end{array}
\right]$
\end{center}

\paragraph{The wedge product.} Consider again M = $\{\varnothing, \{1\}, \{2\}, \{3\}, \{1,2\}, \{1,3\}, \{2,3\}, \{1,2,3\}\}$. Then

$$\varnothing \wedge x = 0 \ \ \forall x$$

$$\{1\} \wedge \{1\} = \{1\}$$
$$\{2\} \wedge \{2\} = \{2\}$$
$$\{3\} \wedge \{3\} = \{3\}$$

$$\{1\} \wedge \{2\} = \{2\} \wedge \{1\} = \{3\} \wedge \{1\} = \{1\} \wedge \{3\} = \{3\} \wedge \{2\} = \{2\} \wedge \{3\} = 0$$

$$\{1\} \wedge \{1,2\} = \{1,2\} \wedge \{1\} = \{2\} \wedge \{1,2\} = \{1,2\} \wedge \{2\} = 1/2 \cdot \{1,2\}$$
$$\{1\} \wedge \{1,3\} = \{1,3\} \wedge \{1\} = \{3\} \wedge \{1,3\} = \{1,3\} \wedge \{3\} = 1/2 \cdot \{1,3\}$$
$$\{2\} \wedge \{2,3\} = \{2,3\} \wedge \{2\} = \{3\} \wedge \{2,3\} = \{2,3\} \wedge \{3\} = 1/2 \cdot \{2,3\}$$

$$\{1\} \wedge \{2,3\} = \{2,3\} \wedge \{1\} = \{2\} \wedge \{1,3\} = \{1,3\} \wedge \{2\} = \{3\} \wedge \{1,2\} = \{1,2\} \wedge \{3\} = 0$$
The first non-trivial products are

$$\{1,2\} \wedge \{1,3\} = - \{1,3\} \wedge \{1,2\} = 1/6 \cdot \{1,2,3\}$$
$$\{1,2\} \wedge \{2,3\} = - \{2,3\} \wedge \{1,2\} = 1/6 \cdot \{1,2,3\}$$
$$\{1,3\} \wedge \{2,3\} = - \{2,3\} \wedge \{1,3\} = 1/6 \cdot \{1,2,3\}$$

\paragraph{Local envelope and localisation.} Consider the same M = $\{\varnothing, \{1\}, \{2\}, \{3\}, \{1,2\}, \{1,3\}, \{2,3\}, \{1,2,3\}\}$. The space of pairs $\Omega_* \otimes \Omega_*$ is spanned by all pairs of simplexes.

Let us provide several examples of the local envelope, defined in s.4.4:
$$\cup (\varnothing, x) = x \ \ \forall x$$

$$\cup (\{1\}, \{1\}) = \{1\}$$

$$\cup (\{1\}, \{2\}) = \{1,2\}$$

$$\cup (\{1\}, \{1,2\}) = \{1,2\}$$

$$\cup (\{1\}, \{2,3\}) = \{1,2,3\}$$

$$\cup (\{1,2\}, \{1,3\}) = \{1,2,3\}$$
and so on. Operations $d$ and $\partial$ act on pairs via lifting (s.5.2):
$$d (a \otimes b) = - d(a) \otimes b - (-1)^{|a|} a \otimes d(b)$$

$$\partial (a \otimes b) = - \partial(a) \otimes b - (-1)^{|a|} a \otimes \partial(b)$$
For example,
$$d (\{1\} \otimes \{1\}) = + \{1,2\} \otimes \{1\} + \{1,3\} \otimes \{1\} + \{1\} \otimes \{1,2\} + \{1\} \otimes \{1,3\}$$

$$d (\{1,2\} \otimes \{1\}) = - \{1,2,3\} \otimes \{1\} - \{1,2\} \otimes \{1,2\} - \{1,2\} \otimes \{1,3\}$$

$$d (\{1\} \otimes \{2,3\}) = \{1,2\} \otimes \{2,3\} + \{1,3\} \otimes \{2,3\} - \{1\} \otimes \{1,2,3\}$$

$$\partial (\{1\} \otimes \{1\}) = - \varnothing \otimes \{1\} - \{1\} \otimes \varnothing$$

$$\partial (\{1,2\} \otimes \{1\}) = \{1\} \otimes \{1\} - \{2\} \otimes \{1\} +  \{1,2\} \otimes \varnothing$$

$$\partial (\{1\} \otimes \{2,3\}) = - \varnothing \otimes \{2,3\} + \{1\} \otimes \{2\} - \{1\} \otimes \{3\}$$
Localisation means omitting the terms of the right hand side, which have local envelope not equal to the local envelope of argument:
$$[d] (\{1\} \otimes \{1\}) = 0$$

$$[d] (\{1,2\} \otimes \{1\}) = - \{1,2\} \otimes \{1,2\}$$

$$[d] (\{1\} \otimes \{2,3\}) = \{1,2\} \otimes \{2,3\} + \{1,3\} \otimes \{2,3\} - \{1\} \otimes \{1,2,3\}$$

$$[\partial] (\{1\} \otimes \{1\}) = - \varnothing \otimes \{1\} - \{1\} \otimes \varnothing$$

$$[\partial] (\{1,2\} \otimes \{1\}) = - \{2\} \otimes \{1\} + \{1,2\} \otimes \varnothing$$

$$[\partial] (\{1\} \otimes \{2,3\}) = 0$$

Everything depends on the number of free vertices in a pair $e = \sigma_1 \otimes \sigma_2$. If there are none, then $[\partial] e = \partial e$. If all vertices in $e$ are free, then $[\partial] e = 0$. Local differential $[d] e = 0$ if there are no free vertices.

Note that there is nothing special in the space of pairs: in the higher spaces $\Omega_*^{\otimes p}$, everything is similar.

\paragraph{Local laplace operator} is equal to $[d] [\partial] + [\partial] [d]$, by definition. For example,

$$\bigtriangleup_{loc} \big(\{1,2\} \otimes \{1\}\big) = [d] \big([\partial](\{1,2\} \otimes \{1\})\big) + [\partial]\big([d](\{1,2\} \otimes \{1\})\big) = $$
$$ = - [d] \big( \{2\} \otimes \{1\} - \{1,2\} \otimes \varnothing \big) - [\partial] \big( \{1,2\} \otimes \{1,2\} \big) = $$
$$ \big( \{1,2\} \otimes \{1\} - \{2\} \otimes \{1,2\} + \{1,2\} \otimes \{1\} + \{1,2\} \otimes \{2\} \big) - \big( \{1\} \otimes \{1,2\} - \{2\} \otimes \{1,2\} - \{1,2\} \otimes \{1\} + \{1,2\} \otimes \{2\} \big) = $$
$$ = 3 \cdot \{1,2\} \otimes \{1\} - \{1\} \otimes \{1,2\}$$

Such calculations lead to the results, described in s.7.3 - 7.4. For example,

\paragraph{1.} On the subspace $\Omega_*^{\otimes 2}(\{1,2\},\{2\})$ with the basis

$$(1,0) = \{1\} \otimes \{1,2\}$$
$$(0,1) = \{1,2\} \otimes \{1\}$$

the local Laplace matrix has a form

\[
\left[ \begin {array}{cccc}
3 & -1 \\
-1 & 3 \\
\end {array} \right]
\]

\begin{center}
eigenvalue 2: eigenvector (1, 1)\end{center}

\begin{center}
eigenvalue 4: eigenvector (-1, 1)\end{center}

\paragraph{2.} On the subspace $\Omega_*^{\otimes 2}(\{1,2,3\},\{2\})$ with the basis

$$(1,0) = \{1,3\} \otimes \{1,2,3\}$$
$$(0,1) = \{1,2,3\} \otimes \{1,3\}$$

the local Laplace matrix has a form

\[
\left[ \begin {array}{cccc}
5 & 1 \\
1 & 5 \\
\end {array} \right]
\]

\begin{center}
eigenvalue 4: eigenvector (-1, 1)\end{center}

\begin{center}
eigenvalue 6: eigenvector (1, 1)\end{center}

\paragraph{3.} On the subspace $\Omega_*^{\otimes 2}(\{1,2\},\{1,2\})$ with the basis

$$(1,0,0,0) = \varnothing \otimes \{1,2\}$$
$$(0,1,0,0) = \{1\} \otimes \{2\}$$
$$(0,0,1,0) = \{2\} \otimes \{1\}$$
$$(0,0,0,1) = \{1,2\} \otimes \varnothing$$

the local Laplace matrix has a form

\[
\left[ \begin {array}{cccc} 2&1&-1&0\\\noalign{\medskip}1&2&0&1\\\noalign{\medskip}-1&0&2&-1\\\noalign{\medskip}0&1&-1&2\end {array} \right]
\]

\begin{center}
eigenvalue 0: eigenvector $(1, -1, 1, 1)$\end{center}

\begin{center}
eigenvalue 2: eigenvectors $(0, 1, 1, 0) \oplus (-1, 0, 0, 1)$\end{center}

\begin{center}
eigenvalue 4: eigenvector $(1, 1, -1, 1)$\end{center}
Construction (\ref{Eigenv}) explains all these results. Compare cases 1 and 2 to observe the non-trivial $\chi$ sign.

\subsection{Discrete quantum field theory}

There is a straightforward physical interpretation of
the construction.
The local Laplace operator $\bigtriangleup_{loc}$ is a
Hamiltonian operator for a discrete quantum field theory,
defined on a Fock phase space
$$\Omega^*_* = \bigoplus_{p = 0}^{\infty}\Omega_*^{\otimes p}$$
which is constructed from a one-particle phase space $\Omega_*$.
The one-particle states are described by forms on a
simplicial complex (wave functions).
The free points correspond to excitations.
Note that the Hamiltonian $\bigtriangleup_{loc}$
arises naturally from the local differentials on complex.
It is essentially discrete and vanishes in the continuum limit.
It would be interesting to investigate the field-theoretical
properties of this theory, compared to smooth analogues.

\subsection{Inverse of Local Laplace Operator}

Inversion of local Laplace operator is necessary to construct the solution (\ref{ZOL2}). We present a combinatorial rule of action of inverse operator $\bigtriangleup_{loc}^{-1}$ on any p-uple of simplexes $e$:

$$\bigtriangleup_{loc}^{-1} e = \Big( [d] [\partial] + [\partial] [d] \Big)^{-1} \Big(\sigma_{1}, \ldots, \sigma_{p}\Big)$$

The result is equal to a sum over all rearrangements of free points into $p$ groups (simplexes). Non-free points remain on their places. Each rearrangement is taken with a weight

\begin{center}
$F(k,n,j = k) = \big(Q^{0}_{0} \big)/F_{0}$
\end{center}
\begin{center}
$F(k,n,j = k-1) = \big(Q^{0}_{1} + Q^{1}_{1} (n-k)\big)/F_{0}$                     \end{center}
\begin{center}
$F(k,n,j = k-2) = \big(Q^{0}_{2} + Q^{1}_{2} (n-k) + Q^{2}_{2} (n-k)(n-k+1)\big)/F_{0}$
\end{center}
\begin{center}
$F(k,n,j = k-3) = \big(Q^{0}_{3} + Q^{1}_{3}(n-k) + Q^{2}_{3} (n-k)(n-k+1) + Q^{3}_{3} (n-k)(n-k+1)(n-k+2)\big)/F_{0}$
\end{center}
\begin{center}
...
\end{center}
where $F_{0} = p^{k} n(n-1)(n-2)...(n-k)$, $Q^{a}_{b} = C^{a}_{b} p^{a} (k-a)!$ and $j$ is the number of free points in the rearrangement that has been moved to another simplex. Also, each rearrangement should be taken with a sign $(-1)^{s1} (-1)^{s2}$. $(-1)^{s1}$ is the parity of permutation of free points, associated with the rearrangement. $(-1)^{s2}$ counts the parity of jumps of free points over the non-free points. This rule certainly provides the same results as the main spectral formula (\ref{Spectral}), but it is more convenient for calculations, for reasons pointed out in s.7.6.

\subsection{Discrete forms and operations: co-chain description}

In the co-chain picture a discrete $p$-form is a function on $p$-simplexes, which can be represented as function $f(i,j,\ldots)$ of $p + 1$ vertexes $i,j, \ldots$ Note that formulas below do not depend on the simplicial complex at all. This is a nice property of local operations in the co-chain description.

\paragraph{The differential.} The differential acts as follows:

\begin{center}
on 0-forms $f(i)$: \ \ $\big(d f\big) (i,j) = f(j) - f(i)$
\end{center}

\begin{center}
on 1-forms $\psi(i,j)$: \ \ $\big(d \psi\big) (i,j,k) = \psi(i,j) - \psi(i,k) + \psi(j,k)$
\end{center}
\begin{center}
on 2-forms $\rho(i,j,k)$: \ \ $\big(d \rho\big) (i,j,k,l) = - \rho(i,j,k) + \rho(i,j,l) - \rho(i,k,l) + \rho(j,k,l)$
\end{center}
and so on.
\paragraph{The wedge product.} The wedge product acts as follows:

\begin{center}
a 0-form $f$ with 0-form $g$: \ \ $\wedge (f, g) (i) = f(i) \cdot g(i)$
\end{center}

\begin{center}
a 0-form $f$ with 1-form $\psi$: \ \ $\wedge (f, \psi) (i,j) = 1/2 \big(f(i) + f(j)\big) \cdot \psi(i,j)$
\end{center}

\begin{center}
a 0-form $f$ with 2-form $\rho$: \ \ $\wedge (f, \rho) (i,j,k) = 1/3 \big(f(i) + f(j) + f(k)\big) \cdot \rho(i,j,k)$
\end{center}

\begin{center}
a 1-form $a$ with 1-form $b$: \ \ $\wedge (a, b) (i,j,k) = 1/6 \big(a(i, k) b(j, k)+a(i, j) b(j, k)-a(j, k) b(i, k)+a(i, j) b(i, k)-a(j, k) b(i, j)-a(i, k) b(i, j) )$
\end{center}

\begin{center}
a 1-form $a$ with 2-form $b$: \ \ $\wedge (a, b) (i,j,k,l) = 1/12 \big(a(i, l) b(j, k, l)+a(i, k) b(j, k, l)+a(i, j) b(j, k, l)-a(j, l) b(i, k, l)-a(j, k) b(i, k, l)+a(i, j) b(i, k, l)+a(k, l) b(i, j, l)-a(j, k) b(i, j, l)-a(i, k) b(i, j, l)+a(k, l) b(i, j, k)+a(j, l) b(i, j, k)+a(i, l) b(i, j, k))$
\end{center}

\begin{center}
a 2-form $a$ with 2-form $b$: \ \ $\wedge (a, b) (i,j,k,l,p) = 1/30 \cdot \big(-a(i, j, l) b(j, k, p)+a(i, j, p) b(i, k, l)+a(k, l, p) b(i, j, p)+a(j, k, l) b(i, j, p)+a(i, k, l) b(i, j, p)+a(k, l, p) b(i, j, l)+a(k, l, p) b(i, j, k)+a(j, l, p) b(i, j, k)+a(j, k, l) b(i, l, p)+a(i, j, k) b(i, l, p)+a(j, k, l) b(i, k, p)+a(i, j, k) b(k, l, p)+a(i, j, k) b(j, l, p)+a(i, l, p) b(j, k, p)+a(i, l, p) b(j, k, l)+a(i, k, p) b(j, k, l)+a(j, k, p) b(i, l, p)+a(i, j, p) b(k, l, p)+a(i, j, l) b(k, l, p)-a(i, k, p) b(j, l, p)-a(j, k, p) b(i, k, l)-a(j, l, p) b(i, k, l)-a(i, j, l) b(i, k, p)-a(j, l, p) b(i, k, p)-a(i, k, l) b(j, k, p)+a(i, l, p) b(i, j, k)+a(i, j, p) b(j, k, l)-a(i, k, l) b(j, l, p)-a(i, k, p) b(i, j, l)-a(j, k, p) b(i, j, l))$
\end{center}
and so on.

\paragraph{The associator.}

Associator of $\wedge$ is defined as $\wedge^2 (a,b,c) = \wedge(a,\wedge(b,c)) - \wedge(\wedge (a,b),c)$. It is a local operation, too. In the co-chain description it can be represented by the following formulas $\big($ three numbers $(n_1,n_2,n_3)$ denote the degrees of forms $a,b,c$ respectively $\big)$

\begin{center}
(1,0,0): \ \ $\wedge^2 (a, b, c) (i,j) = 1/4 \cdot a(i, j) b(j) c(j) + 1/4 \cdot a(i, j) b(i) c(i) - 1/4 \cdot c(j) a(i, j) b(i) - 1/4 \cdot c(i) a(i, j) b(j)$
\end{center}

\begin{center}
(0,1,1): \ \ $\wedge^2 (a, b, c) (i,j,k) = -1/36 \cdot a(k) b(i, k) c(j, k)+1/18 \cdot a(k) b(i, j) c(j, k)+1/36 \cdot a(k) b(j, k) c(i, k)+1/18 \cdot a(k) b(i, j) c(i, k)+1/36 \cdot a(k) b(j, k) c(i, j)+1/18 \cdot a(j) b(i, k) c(j, k)-1/36 \cdot a(j) b(i, j) c(j, k)+1/36 \cdot a(j) b(j, k) c(i, k)-1/36 \cdot a(j) b(i, j) c(i, k)+1/36 \cdot a(j) b(j, k) c(i, j)-1/18 \cdot a(j) b(i, k) c(i, j)-1/36 \cdot a(i) b(i, k) c(j, k)-1/36 \cdot a(i) b(i, j) c(j, k)-1/18 \cdot a(i) b(j, k) c(i, k)-1/36 \cdot a(i) b(i, j) c(i, k)-1/18 \cdot a(i) b(j, k) c(i, j)+1/36 \cdot a(i) b(i, k) c(i, j)+1/36 \cdot a(k) b(i, k) c(i, j)$
\end{center}

\begin{center}
(1,1,1): \ \ $\wedge^2 (a, b, c) (i,j,k,l) = -1/72 \cdot c(i, j) a(i, k) b(k, l)+1/72 \cdot c(k, l) a(i, k) b(i, j)-1/72 \cdot c(j, l) a(k, l) b(i, k)-1/72 \cdot c(j, l) a(i, l) b(i, k)+1/72 \cdot c(j, l) a(j, k) b(i, k)-1/72 \cdot c(j, l) a(i, j) b(i, k)-1/72 \cdot c(j, k) a(k, l) b(i, l)+1/72 \cdot c(j, k) a(i, k) b(i, l)-1/72 \cdot c(j, k) a(j, l) b(i, l)+1/72 \cdot c(j, k) a(i, j) b(i, l)+1/72 \cdot c(i, l) a(k, l) b(j, k)+1/72 \cdot c(i, l) a(j, l) b(j, k)-1/72 \cdot c(i, l) a(i, k) b(j, k)-1/72 \cdot c(i, l) a(i, j) b(j, k)+1/72 \cdot c(i, k) a(k, l) b(j, l)-1/72 \cdot c(i, k) a(j, k) b(j, l)+1/72 \cdot c(i, k) a(i, l) b(j, l)+1/72 \cdot c(i, k) a(i, j) b(j, l)-1/72 \cdot c(i, j) a(j, l) b(k, l)-1/72 \cdot c(i, j) a(j, k) b(k, l)-1/72 \cdot c(i, j) a(i, l) b(k, l)+1/72 \cdot a(i, l) b(j, k) c(k, l)+1/72 \cdot a(i, l) b(j, k) c(j, l)+1/72 \cdot a(i, k) b(j, l) c(k, l)-1/72 \cdot a(i, k) b(j, l) c(j, k)-1/72 \cdot a(i, j) b(k, l) c(j, l)-1/72 \cdot a(i, j) b(k, l) c(j, k)-1/72 \cdot a(j, l) b(i, k) c(k, l)-1/72 \cdot a(j, l) b(i, k) c(i, l)-1/72 \cdot a(j, k) b(i, l) c(k, l)+1/72 \cdot a(j, k) b(i, l) c(i, k)-1/72 \cdot a(i, j) b(k, l) c(i, l)-1/72 \cdot a(i, j) b(k, l) c(i, k)+1/72 \cdot a(k, l) b(i, j) c(j, l)+1/72 \cdot a(k, l) b(i, j) c(i, l)-1/72 \cdot a(j, k) b(i, l) c(j, l)+1/72 \cdot a(j, k) b(i, l) c(i, j)+1/72 \cdot a(i, k) b(j, l) c(i, l)+1/72 \cdot a(i, k) b(j, l) c(i, j)+1/72 \cdot a(k, l) b(i, j) c(j, k)+1/72 \cdot a(k, l) b(i, j) c(i, k)+1/72 \cdot a(j, l) b(i, k) c(j, k)-1/72 \cdot a(j, l) b(i, k) c(i, j)-1/72 \cdot a(i, l) b(j, k) c(i, k)-1/72 \cdot a(i, l) b(j, k) c(i, j)+1/72 \cdot c(k, l) a(j, l) b(i, j)+1/72 \cdot c(k, l) a(i, l) b(i, j)+1/72 \cdot c(k, l) a(j, k) b(i, j)$
\end{center}
and so on. Associator has definite symmetry properties, that follow from the skew-symmetry of $\wedge$ product.

\subsection{Higher operation $m = m^{(3)}$}

Higher operation $m = m^{(3)}$ is defined as a solution to

$$dm(a,b,c) + m(da,b,c) + (-1)^{|a|} m(a,db,c) +
(-1)^{|a| + |b|} m(a,b,dc) = \big{(}(a \wedge b) \wedge c - a \wedge (b \wedge c)\big{)}$$
This paper is devoted to solving such equations. First of all, let us consider an example of naive solution obtained with the help of (\ref{SOL1}). Naive solution is not local and depends on the simplicial complex. For example, consider a 2-disc $M = \{\varnothing,\{1\}, \{2\}, \{3\}, \{1,2\}, \{1,3\}, \{2,3\}, \{1,2,3\}\}$. The naive solution is $m (a, b, c)_{\sigma} = - 1/3 \cdot \partial \wedge^2 (a,b,c)_{\sigma}$. Explicitly, if $\phi$ is a 0-form and $\psi_1, \psi_2$ are 1-forms, then

$$m (\phi, \psi_1, \psi_2) (1,2) = 1/3 \cdot \wedge^2 (f, \psi_1, \psi_2) (1,2,3)$$
$$m (\phi, \psi_1, \psi_2) (1,3) = - 1/3 \cdot \wedge^2 (f, \psi_1, \psi_2) (1,2,3)$$
$$m (\phi, \psi_1, \psi_2) (2,3) = 1/3 \cdot \wedge^2 (f, \psi_1, \psi_2) (1,2,3)$$

To calculate the value of $m(\phi, \psi_1, \psi_2)$ on simplex $\{1,2\}$ we need to know the values of forms $\phi, \psi_1, \psi_2$ on simplexes beyond $\{1,2\}$, actually, we need to know their values on the whole complex. This is another illustration of non-locality.

Second, we present the local solution $m^{(3)}$, obtained with the help of (\ref{ZOL2}). The numeric factors (like $7/1296$) typically arise from the inversion of local Laplace operator.

\begin{center}
(0,1,1): \ \ $m (a, b, c) (i,j) = 1/8 \cdot a(j) b(i,j) c(i,j) - 1/8 \cdot a(i) b(i,j) c(i,j)$
\end{center}

\begin{center}
(1,1,1): \ \ $m (a, b, c) (i,j,k) = -1/54 \cdot a(j,k) b(j,k) c(i,k)+1/54 \cdot a(j,k) b(j,k) c(i,j)+1/27 \cdot a(j,k) b(i,k) c(j,k)+1/54 \cdot a(j,k) b(i,k) c(i,k)-1/27 \cdot a(j,k) b(i,j) c(j,k)+1/54 \cdot a(j,k) b(i,j) c(i,j)-1/54 \cdot a(i,k) b(j,k) c(j,k)-1/27 \cdot a(i,k) b(j,k) c(i,k)+1/54 \cdot a(i,k) b(i,k) c(j,k)+1/54 \cdot a(i,k) b(i,k) c(i,j)-1/27 \cdot a(i,k) b(i,j) c(i,k)-1/54 \cdot a(i,k) b(i,j) c(i,j)+1/54 \cdot a(i,j) b(j,k) c(j,k)-1/27 \cdot a(i,j) b(j,k) c(i,j)+1/54 \cdot a(i,j) b(i,k) c(i,k)+1/27 \cdot a(i,j) b(i,k) c(i,j)+1/54 \cdot a(i,j) b(i,j) c(j,k)-1/54 \cdot a(i,j) b(i,j) c(i,k)$
\end{center}

\begin{center}
(2,1,1): \ \ $m (a, b, c) (i,j,k,l) = 1/864 \cdot a(i,j,l) b(j,k) c(i,k)+1/864 \cdot a(j,k,l) b(i,l) c(i,j)+1/864 \cdot a(i,k,l) b(i,j) c(j,l)+17/1296 \cdot a(j,k,l) b(i,l) c(k,l)+1/864 \cdot a(i,k,l) b(j,l) c(i,j)+23/2592 \cdot a(i,k,l) b(k,l) c(j,l)-1/864 \cdot a(i,k,l) b(j,l) c(j,k)+1/864 \cdot a(j,k,l) b(i,l) c(i,k)+1/864 \cdot a(j,k,l) b(i,k) c(i,j)-1/864 \cdot a(i,j,l) b(k,l) c(j,k)-17/1296 \cdot a(i,k,l) b(i,j) c(i,l)-1/864 \cdot a(i,j,l) b(j,k) c(k,l)-17/1296 \cdot a(i,k,l) b(j,l) c(k,l)-17/1296 \cdot a(j,k,l) b(i,j) c(j,l)+7/1296 \cdot a(i,j,l) b(i,k) c(i,k)+1/864 \cdot a(i,k,l) b(j,k) c(i,j)+1/864 \cdot a(i,j,l) b(i,k) c(j,k)-23/2592 \cdot a(i,j,k) b(j,k) c(k,l)+17/1296 \cdot a(i,j,l) b(i,k) c(i,l)-17/1296 \cdot a(i,j,k) b(i,l) c(i,k)+17/1296 \cdot a(i,k,l) b(j,k) c(k,l)-17/1296 \cdot a(i,j,k) b(j,l) c(j,k)+23/2592 \cdot a(j,k,l) b(j,k) c(i,j)+17/1296 \cdot a(j,k,l) b(i,k) c(j,k)-17/1296 \cdot a(i,j,l) b(j,k) c(i,j)+17/1296 \cdot a(i,j,l) b(k,l) c(i,l)-7/1296 \cdot a(i,j,k) b(i,l) c(i,l)-23/2592 \cdot a(i,j,l) b(i,l) c(i,k)+17/1296 \cdot a(i,j,k) b(j,l) c(i,j)+23/2592 \cdot a(i,k,l) b(i,k) c(j,k)-23/2592 \cdot a(i,j,l) b(i,j) c(i,k)-7/1296 \cdot a(i,k,l) b(j,l) c(j,l)-23/2592 \cdot a(j,k,l) b(k,l) c(i,l)-1/864 \cdot a(i,j,k) b(j,l) c(k,l)-7/1296 \cdot a(i,k,l) b(j,k) c(j,k)-23/2592 \cdot a(i,j,k) b(i,j) c(j,l)+1/864 \cdot a(j,k,l) b(i,j) c(i,k)-23/2592 \cdot a(j,k,l) b(j,k) c(i,k)-17/1296 \cdot a(i,k,l) b(i,j) c(i,k)+17/1296 \cdot a(j,k,l) b(i,l) c(j,l)-23/2592 \cdot a(i,j,k) b(i,k) c(k,l)-1/864 \cdot a(i,j,k) b(i,l) c(j,l)-23/2592 \cdot a(i,j,l) b(j,l) c(j,k)+1/864 \cdot a(i,k,l) b(i,j) c(j,k)+7/1296 \cdot a(i,j,l) b(k,l) c(k,l)+17/1296 \cdot a(i,j,k) b(k,l) c(i,k)-7/1296 \cdot a(i,k,l) b(i,j) c(i,j)+17/1296 \cdot a(i,j,k) b(k,l) c(j,k)-7/1296 \cdot a(i,j,k) b(j,l) c(j,l)+1/864 \cdot a(j,k,l) b(i,k) c(i,l)-17/1296 \cdot a(j,k,l) b(i,k) c(k,l)-23/2592 \cdot a(i,j,l) b(j,l) c(k,l)+23/2592 \cdot a(j,k,l) b(j,l) c(i,j)+17/1296 \cdot a(i,j,l) b(k,l) c(j,l)+23/2592 \cdot a(i,j,k) b(i,k) c(i,l)+1/864 \cdot a(j,k,l) b(i,j) c(i,l)-23/2592 \cdot a(i,k,l) b(k,l) c(j,k)-1/864 \cdot a(i,j,k) b(i,l) c(k,l)+7/1296 \cdot a(i,j,l) b(j,k) c(j,k)-1/864 \cdot a(i,j,l) b(k,l) c(i,k)+23/2592 \cdot a(i,j,k) b(i,j) c(i,l)-17/1296 \cdot a(i,k,l) b(j,k) c(i,k)-23/2592 \cdot a(j,k,l) b(j,l) c(i,l)+7/1296 \cdot a(j,k,l) b(i,j) c(i,j)+23/2592 \cdot a(i,k,l) b(i,l) c(j,l)+23/2592 \cdot a(i,k,l) b(i,k) c(i,j)+7/1296 \cdot a(j,k,l) b(i,l) c(i,l)+23/2592 \cdot a(i,k,l) b(i,l) c(i,j)+17/1296 \cdot a(i,j,l) b(i,k) c(i,j)-23/2592 \cdot a(i,j,l) b(i,l) c(k,l)-17/1296 \cdot a(j,k,l) b(i,j) c(j,k)-1/864 \cdot a(i,j,k) b(k,l) c(j,l)-1/864 \cdot a(i,j,l) b(i,k) c(k,l)-7/1296 \cdot a(i,j,k) b(k,l) c(k,l)+23/2592 \cdot a(i,j,k) b(j,k) c(j,l)+23/2592 \cdot a(j,k,l) b(k,l) c(i,k)+7/1296 \cdot a(j,k,l) b(i,k) c(i,k)+23/2592 \cdot a(i,j,l) b(i,j) c(j,k)+17/1296 \cdot a(i,j,l) b(j,k) c(j,l)-1/864 \cdot a(i,j,k) b(j,l) c(i,l)-1/864 \cdot a(i,k,l) b(j,k) c(j,l)-1/864 \cdot a(i,j,k) b(k,l) c(i,l)-17/1296 \cdot a(i,j,k) b(i,l) c(i,j)-17/1296 \cdot a(i,k,l) b(j,l) c(i,l)$
\end{center}

\begin{center}
(1,2,1): \ \ $m (a, b, c) (i,j,k,l) = -1/432 \cdot a(k,l) b(i,j,l) c(i,k)+11/2592 \cdot a(j,l) b(i,j,k) c(i,j)-1/432 \cdot a(i,k) b(i,j,l) c(k,l)-11/2592 \cdot a(i,j) b(i,k,l) c(i,k)-1/432 \cdot a(j,k) b(i,k,l) c(j,l)+1/432 \cdot a(j,k) b(i,k,l) c(i,j)+1/432 \cdot a(j,k) b(i,j,l) c(i,k)+1/432 \cdot a(i,j) b(j,k,l) c(i,l)+11/2592 \cdot a(j,l) b(j,k,l) c(i,l)+1/432 \cdot a(i,k) b(j,k,l) c(i,j)-11/2592 \cdot a(j,k) b(i,j,l) c(i,j)-11/2592 \cdot a(i,j) b(i,k,l) c(i,l)-1/432 \cdot a(j,l) b(i,k,l) c(j,k)+1/432 \cdot a(i,j) b(i,k,l) c(j,l)+1/432 \cdot a(i,j) b(i,k,l) c(j,k)-11/2592 \cdot a(i,j) b(i,j,l) c(j,k)+7/648 \cdot a(k,l) b(i,j,l) c(k,l)-11/2592 \cdot a(i,j) b(j,k,l) c(j,l)-1/432 \cdot a(k,l) b(i,j,k) c(i,l)-1/432 \cdot a(i,l) b(i,j,k) c(k,l)-11/2592 \cdot a(i,k) b(i,k,l) c(j,k)+11/2592 \cdot a(i,l) b(j,k,l) c(k,l)+11/2592 \cdot a(k,l) b(j,k,l) c(i,l)+11/2592 \cdot a(j,k) b(i,k,l) c(k,l)+11/2592 \cdot a(i,k) b(j,k,l) c(j,k)+11/2592 \cdot a(i,k) b(i,j,l) c(i,j)+11/2592 \cdot a(i,j) b(i,j,k) c(j,l)-11/2592 \cdot a(i,k) b(i,k,l) c(i,j)+1/432 \cdot a(j,l) b(i,k,l) c(i,j)-11/2592 \cdot a(i,l) b(i,j,k) c(i,k)+7/648 \cdot a(i,l) b(j,k,l) c(i,l)-11/2592 \cdot a(i,k) b(j,k,l) c(k,l)+11/2592 \cdot a(k,l) b(i,j,k) c(i,k)+11/2592 \cdot a(i,k) b(i,j,k) c(k,l)-11/2592 \cdot a(j,l) b(i,j,k) c(j,k)-11/2592 \cdot a(i,l) b(i,k,l) c(i,j)-1/432 \cdot a(j,l) b(i,j,k) c(i,l)-7/648 \cdot a(i,j) b(i,k,l) c(i,j)-11/2592 \cdot a(j,k) b(j,k,l) c(i,j)+1/432 \cdot a(i,k) b(i,j,l) c(j,k)+11/2592 \cdot a(j,k) b(i,j,k) c(k,l)-11/2592 \cdot a(i,l) b(i,j,k) c(i,j)-11/2592 \cdot a(j,l) b(i,k,l) c(i,l)+1/432 \cdot a(i,j) b(j,k,l) c(i,k)+11/2592 \cdot a(i,k) b(i,j,l) c(i,l)+11/2592 \cdot a(k,l) b(i,k,l) c(j,k)+11/2592 \cdot a(i,l) b(j,k,l) c(j,l)-11/2592 \cdot a(j,l) b(j,k,l) c(i,j)+7/648 \cdot a(i,j) b(j,k,l) c(i,j)-7/648 \cdot a(j,k) b(i,k,l) c(j,k)-7/648 \cdot a(i,l) b(i,j,k) c(i,l)-1/432 \cdot a(k,l) b(i,j,k) c(j,l)+7/648 \cdot a(i,k) b(i,j,l) c(i,k)-1/432 \cdot a(k,l) b(i,j,l) c(j,k)-7/648 \cdot a(j,l) b(i,j,k) c(j,l)+11/2592 \cdot a(k,l) b(i,j,l) c(i,l)-11/2592 \cdot a(j,k) b(i,k,l) c(i,k)+11/2592 \cdot a(j,l) b(i,j,l) c(k,l)+1/432 \cdot a(i,k) b(j,k,l) c(i,l)+7/648 \cdot a(i,k) b(j,k,l) c(i,k)+11/2592 \cdot a(j,l) b(i,j,l) c(j,k)-7/648 \cdot a(k,l) b(i,j,k) c(k,l)-1/432 \cdot a(j,k) b(i,j,l) c(k,l)-11/2592 \cdot a(j,k) b(i,j,k) c(j,l)-11/2592 \cdot a(k,l) b(j,k,l) c(i,k)-11/2592 \cdot a(k,l) b(i,k,l) c(j,l)-11/2592 \cdot a(i,k) b(i,j,k) c(i,l)+11/2592 \cdot a(k,l) b(i,j,l) c(j,l)-11/2592 \cdot a(i,l) b(i,k,l) c(j,l)-1/432 \cdot a(j,l) b(i,j,k) c(k,l)+7/648 \cdot a(j,k) b(i,j,l) c(j,k)-11/2592 \cdot a(i,j) b(j,k,l) c(j,k)-11/2592 \cdot a(i,j) b(i,j,k) c(i,l)+11/2592 \cdot a(j,k) b(j,k,l) c(i,k)+1/432 \cdot a(i,l) b(j,k,l) c(i,k)+11/2592 \cdot a(k,l) b(i,j,k) c(j,k)+11/2592 \cdot a(i,l) b(i,j,l) c(k,l)-1/432 \cdot a(i,l) b(i,j,k) c(j,l)+11/2592 \cdot a(i,j) b(i,j,l) c(i,k)-11/2592 \cdot a(j,l) b(i,k,l) c(k,l)+11/2592 \cdot a(i,l) b(i,j,l) c(i,k)+11/2592 \cdot a(j,k) b(i,j,l) c(j,l)-7/648 \cdot a(j,l) b(i,k,l) c(j,l)+1/432 \cdot a(i,l) b(j,k,l) c(i,j)
$
\end{center}

\begin{center}
(1,1,2): \ \ $m (a, b, c) (i,j,k,l) = 1/864 \cdot a(j,k) b(i,j) c(i,k,l)+1/864 \cdot a(j,l) b(i,j) c(i,k,l)-1/864 \cdot a(j,k) b(j,l) c(i,k,l)+17/1296 \cdot a(i,l) b(k,l) c(i,j,l)+1/864 \cdot a(i,k) b(i,j) c(j,k,l)+1/864 \cdot a(i,l) b(i,j) c(j,k,l)-1/864 \cdot a(k,l) b(i,l) c(i,j,k)-1/864 \cdot a(i,l) b(k,l) c(i,j,k)-1/864 \cdot a(k,l) b(j,k) c(i,j,l)+1/864 \cdot a(i,j) b(j,k) c(i,k,l)+1/864 \cdot a(i,j) b(i,l) c(j,k,l)-1/864 \cdot a(j,l) b(k,l) c(i,j,k)+1/864 \cdot a(i,l) b(i,k) c(j,k,l)+7/1296 \cdot a(i,k) b(i,k) c(j,k,l)-7/1296 \cdot a(k,l) b(k,l) c(i,j,k)-23/2592 \cdot a(j,k) b(j,l) c(i,j,l)+7/1296 \cdot a(i,l) b(i,l) c(j,k,l)-17/1296 \cdot a(i,k) b(j,k) c(i,k,l)+17/1296 \cdot a(k,l) b(i,l) c(j,k,l)-17/1296 \cdot a(i,j) b(i,l) c(i,j,k)+23/2592 \cdot a(j,l) b(k,l) c(i,k,l)+17/1296 \cdot a(i,l) b(i,k) c(i,j,l)+23/2592 \cdot a(i,j) b(j,l) c(j,k,l)-7/1296 \cdot a(j,l) b(j,l) c(i,j,k)+7/1296 \cdot a(i,j) b(i,j) c(j,k,l)-17/1296 \cdot a(i,k) b(i,l) c(i,j,k)+23/2592 \cdot a(j,k) b(i,k) c(i,k,l)+17/1296 \cdot a(j,l) b(i,l) c(j,k,l)-23/2592 \cdot a(i,l) b(k,l) c(j,k,l)-17/1296 \cdot a(i,k) b(i,j) c(i,k,l)+17/1296 \cdot a(j,k) b(k,l) c(i,j,k)-23/2592 \cdot a(i,l) b(j,l) c(j,k,l)-17/1296 \cdot a(j,k) b(j,l) c(i,j,k)-1/864 \cdot a(k,l) b(j,l) c(i,j,k)-23/2592 \cdot a(j,l) b(i,j) c(i,j,k)-17/1296 \cdot a(j,l) b(i,j) c(j,k,l)+17/1296 \cdot a(j,k) b(i,k) c(j,k,l)-7/1296 \cdot a(j,l) b(j,l) c(i,k,l)+23/2592 \cdot a(i,j) b(i,l) c(i,k,l)-7/1296 \cdot a(i,l) b(i,l) c(i,j,k)-17/1296 \cdot a(k,l) b(i,k) c(j,k,l)+7/1296 \cdot a(i,k) b(i,k) c(i,j,l)+1/864 \cdot a(i,k) b(i,l) c(j,k,l)-23/2592 \cdot a(k,l) b(j,l) c(i,j,l)+23/2592 \cdot a(i,l) b(i,j) c(i,j,k)+17/1296 \cdot a(i,j) b(i,k) c(i,j,l)-17/1296 \cdot a(k,l) b(j,l) c(i,k,l)-23/2592 \cdot a(k,l) b(i,k) c(i,j,k)-23/2592 \cdot a(i,k) b(i,l) c(i,j,l)+23/2592 \cdot a(i,j) b(i,k) c(i,k,l)-23/2592 \cdot a(i,k) b(j,k) c(j,k,l)-17/1296 \cdot a(i,l) b(j,l) c(i,k,l)+17/1296 \cdot a(j,l) b(k,l) c(i,j,l)-17/1296 \cdot a(j,k) b(i,j) c(j,k,l)-1/864 \cdot a(j,l) b(j,k) c(i,k,l)+1/864 \cdot a(i,j) b(j,l) c(i,k,l)-23/2592 \cdot a(i,k) b(i,j) c(i,j,l)-1/864 \cdot a(i,k) b(k,l) c(i,j,l)+7/1296 \cdot a(k,l) b(k,l) c(i,j,l)+1/864 \cdot a(i,k) b(j,k) c(i,j,l)+17/1296 \cdot a(j,l) b(j,k) c(i,j,l)+17/1296 \cdot a(i,k) b(k,l) c(i,j,k)+23/2592 \cdot a(j,l) b(j,k) c(i,j,k)-1/864 \cdot a(k,l) b(i,k) c(i,j,l)-23/2592 \cdot a(j,k) b(k,l) c(i,k,l)+23/2592 \cdot a(j,k) b(i,j) c(i,j,l)+23/2592 \cdot a(i,j) b(j,k) c(j,k,l)+17/1296 \cdot a(i,j) b(j,l) c(i,j,k)+1/864 \cdot a(j,k) b(i,k) c(i,j,l)+23/2592 \cdot a(j,l) b(i,l) c(i,k,l)-23/2592 \cdot a(k,l) b(j,k) c(i,j,k)+7/1296 \cdot a(j,k) b(j,k) c(i,j,l)+1/864 \cdot a(i,j) b(i,k) c(j,k,l)-23/2592 \cdot a(k,l) b(i,l) c(i,j,l)-17/1296 \cdot a(i,l) b(i,j) c(i,k,l)-7/1296 \cdot a(i,j) b(i,j) c(i,k,l)-7/1296 \cdot a(j,k) b(j,k) c(i,k,l)+17/1296 \cdot a(k,l) b(j,k) c(i,k,l)-1/864 \cdot a(j,k) b(k,l) c(i,j,l)-1/864 \cdot a(i,l) b(j,l) c(i,j,k)-1/864 \cdot a(j,l) b(i,l) c(i,j,k)+23/2592 \cdot a(i,k) b(k,l) c(j,k,l)-17/1296 \cdot a(i,j) b(j,k) c(i,j,l)+23/2592 \cdot a(i,l) b(i,k) c(i,j,k)$
\end{center}
These formulas represent the local operation $m^{(3)}$ - a local solution to $A_\infty$ equations. Higher expressions become more complicated. However, the solution (even if lengthy) is always straightforward, because the explicit formula (\ref{ZOL2}) is written.

\section{Acknowledgements}

Our work is partly supported by Federal Nuclear Energy Agency,
by the joint grant 06-01-92059-CE,  by NWO project
047.011.2004.026, by INTAS grant 05-1000008-7865,
by ANR-05-BLAN-0029-01 project (A.M.), by RFBR grants 07-02-00642 (A.M. and Sh.Sh) and 07-02-01161 (V.D.), by the Grant of Support for the Scientific Schools NSh-8004.2006.2.

\end{document}